\documentclass{amsart}

\usepackage[utf8]{inputenc}
\usepackage[english]{babel}

\usepackage{amssymb}
\usepackage{amsfonts}
\usepackage{amscd}

\newcommand{\R}{{\mathbb{R}}}

\newcommand{\Dbb}{{\mathbb{D}}}
\newcommand{\Ubb}{{\mathbb{U}}}

\newcommand{\Ac}{{\mathcal{A}}}
\newcommand{\Po}{{\mathcal{P}}}
\newcommand{\Fo}{{\mathcal{F}}}

\newcommand{\Hc}{{\mathcal{H}}}

\newcommand{\Gc}{{\mathcal{G}}}

\newcommand{\de}{{\mathrm{d}}}
\newcommand{\vol}{{\mathrm{vol}}}
\newtheorem{theorem}{Theorem}[section]
\newtheorem{corollary}[theorem]{Corollary}

\newtheorem{definition}[theorem]{Definition}

\theoremstyle{definition}
\theoremstyle{remark}
\numberwithin{equation}{section}

\newcommand{\translation}[1]{$[\![$#1$]\!]$}

\newcommand{\labelpag}[1]{\refstepcounter{enumi}\label{#1}}

\begin{document}
\title[Peano on derivative of measures]{
Peano on derivative of measures: \\
strict derivative of 
distributive set functions
}
\author{Gabriele H. Greco}
\address{Dipartimento di Matematica\\
Universit\`{a} di Trento, 38050 Povo (TN), Italy}
\email{greco@science.unitn.it}
\author{Sonia Mazzucchi}
\address{Dipartimento di Matematica\\
Universit\`{a} di Trento, 38050 Povo (TN), Italy}
\email{mazzucch@science.unitn.it}
\author{Enrico M. Pagani}
\address{Dipartimento di Matematica\\
Universit\`{a} di Trento, 38050 Povo (TN), Italy}
\email{pagani@science.unitn.it}

\dedicatory{On the occasion of the $150^{th}$ anniversary of the birth of Giuseppe Peano}
%Commemorating the 150th birthday of Giuseppe Peano (1858-1932)}
\date{%
February 22, 2010% (submitted to Colloquium Mathematicum)%\today
}

\begin{abstract}

By retracing research on \emph{coexistent magnitudes}  (\emph{grandeurs coexistantes}) by \textsc{Cauchy} \cite[(1841)]{cauchy1841}, \textsc{Peano}  in \emph{Applicazioni geometriche del calcolo infinitesimale} \cite[(1887)]{peano87} defines the ``density'' (strict derivative) of a ``mass" (a distributive set function) with respect to a ``volume'' (a positive distributive set function), proves its continuity (whenever the  
strict derivative exists) and shows the validity of the \emph{mass-density paradigm}:  ``mass'' is recovered from  ``density'' by integration with respect to ``volume''.
It is remarkable that \textsc{Peano}'s strict derivative provides a consistent mathematical ground to the concept of ``infinitesimal ratio'' between two magnitudes, successfully used since \textsc{Kepler}.
In this way  the classical (i.e., pre-Lebesgue) measure theory reaches a complete and definitive form in \textsc{Peano}'s
\emph{Applicazioni geometriche}.

A primary aim of the present paper is a detailed exposition of \textsc{Peano}'s work of 1887 leading to the concept of strict derivative of distributive set functions and their use.
Moreover, we compare \textsc{Peano}'s work and \textsc{Lebesgue}'s \emph{La mesure des grandeurs}
\cite[(1935)]{lebesgue1935}: in this memoir \textsc{Lebesgue}, motivated by coexistent magnitudes of \textsc{Cauchy}, introduces a uniform-derivative of certain additive set functions, a concept that coincides with
\textsc{Peano}'s strict derivative. Intriguing questions are whether \textsc{Lebesgue} was aware of the
contributions of \textsc{Peano} and which role is played by the notions of strict derivative or of uniform-derivative in today mathematical practice.

\end{abstract}

\maketitle
\section{Introduction}
By referring to \textsc{Cauchy} \cite[(1841)]{cauchy1841} \textsc{Peano} introduces in \emph{Applicazioni geometriche del calcolo infinitesimale} \cite[(1887)]{peano87} the concept of \emph{strict derivative} of set functions. The set functions 
considered by him are not precisely  finite additive measures.  The modern concept of finite additivity is based on partitions by disjoint sets, while \textsc{Peano}'s additivity property coincides with a traditional
supple concept of ``decompositions of magnitudes'', which \textsc{Peano} implements in his proofs as 
\emph{distributive set functions}.

Contrary to \textsc{Peano}'s strict derivative (\emph{rapporto}), \textsc{Cauchy}'s derivative (\emph{rapport diff\'erentiel}) of a set function corresponds to the usual derivative of functions of one variable. In \textsc{Peano}'s Theorem \ref{peanodev} on strict derivative of distributive set functions the (physical) \emph{mass-density paradigm} is realized: the ``mass'' (a distributive set function) is recovered from the ``density'' (strict derivative) by integration with respect to the ``volume'' (a positive distributive set function of reference).

\textsc{Peano}  expresses  \textsc{Cauchy}'s ideas in a  more precise and modern language and completes the program proposed by \textsc{Cauchy}, who, at the end of his article \cite[(1841) p.\,229]{cauchy1841}, writes: 
\begin{quote}
Dans un autre M\'emoire  nous donnerons de nouveaux d\'eveloppements aux principes ci-dessus expos\'es [on coexistent magnitudes], en les appliquant d'une mani\`ere sp\'eciale \a l'\'evalutation des longueurs, des aires et des volumes.\,\footnote{\translation{%
In another memoir we will give new developments to the above mentioned statements  [on coexistent magnitudes], and we will apply them to evaluate lengths, areas and volumes.}}
\end{quote}

Among numerous applications of \textsc{Peano}'s strict derivatives of set functions which can be found in \emph{Applicazioni geometriche}, there are formulae on oriented integrals, in which the geometric vector calculus by \textsc{Grassmann} plays an important role. For instance, \textsc{Peano} proves the  formula of area 
starting by his definition of area of a surface, that he proposed in order to solve the drawbacks of \textsc{Serret}'s definition of area \cite[(1879)]{serret}.

The didactic value of  \textsc{Peano}'s strict derivative of set function is transparent: in \emph{La mesure
des grandeurs} \cite[(1935)]{lebesgue1935} \textsc{Lebesgue} himself  uses a similar approach to differentiation of measures in order to simplify the exposition of his measure theory.

In  Section \ref{sez-paradigma}, \textsc{Peano}'s and \textsc{Lebesgue}'s derivative are compared in view  of the  paradigm of mass-density  
and of the paradigm of primitives, that motivated  mathematical research between  $19^{th}$ century and the beginning of $20^{th}$ century.
In the celebrated paper \emph{L'int\'egration des fonctions discontinues}  \cite[(1910)]{lebesgue1910}  \textsc{Lebesgue} defines a derivative of  $\sigma$-additive measures with respect to the volume. He proves its existence 
and its measurability. In the case of absolute continuity of  the $\sigma$-additive measures, \textsc{Lebesgue} proves that the measure is given by the integral of his derivative with respect to the volume.
As it will be seen later in details,  \textsc{Peano}'s \emph{strict derivative} of distributive set functions  does not necessarily exist and, moreover, whenever it exists, \textsc{Peano}'s strict derivative is continuous, while \textsc{Lebesgue}'s derivative in general is not.

Section \ref{sez-misura} presents an overview of \textsc{Peano}'s work on pre-Lebesgue classical measure theory which is completed in Sections \ref{sez-distributive}-\ref{sez-derivata}.

Section \ref{sez-cauchy} is devoted to an analysis of \textsc{Cauchy}'s  \emph{Coexistent magnitudes}
\cite[(1841)]{cauchy1841}\,%
\footnote{From now on we refer to \textsc{Cauchy}'s paper \emph{M{\'e}moire sur le rapport diff{\'e}rentiel de deux grandeurs qui
  varient simultan{\'e}ment} \cite[(1841)]{cauchy1841} as to \emph{Coexistent magnitudes.}
}, by emphasizing the results that will be found, in a different language, 
in \textsc{Peano}'s \emph{Applicazioni geometriche}
or in \textsc{Lebesgue}'s \emph{La mesure des grandeurs}. 

Section \ref{sez-distributive} concerns the concept of ``distributive families" and of ``distributive set functions'' as presented by \textsc{Peano}  in \emph{Applicazioni geometriche} and in his paper \emph{Le grandezze coesistenti di Cauchy} \cite[(1915)]{peano1915}. 

Section \ref{sez-derivata} presents a definition of strict derivative of set functions,  main results and some applications, while in  Section \ref{sez-massdensity} we discuss \textsc{Peano}'s 
definition of  integral of set functions and a related theorem that realizes the mentioned
physical paradigm of  mass-density.

Section \ref{sez-comments} presents  the approach of \textsc{Lebesgue} in \emph{La mesure des grandeurs} to  \textsc{Cauchy}'s coexistent magnitudes,  leading to  introduction of a new notion of derivative: the uniform-derivative.

We observe that this paper is meanly historical.
From a methodological point of view, we are focussed on primary sources, that is, on mathematical facts and not on the elaborations or interpretations of these facts by other scholars of history of mathematics. 
For convenience of the reader, original statements and, in some case, terminology are presented
in a modern form, preserving, of course, their content.

Historical investigations on forgotten mathematical achievements are not useless (from the point of view of mathematics),
because some of them carry ideas that remain innovative today. This thought was very well expressed by
 \textsc{Mascheroni}  before the beginning of the study of the geometrical problems leading to the \emph{Geo\-me\-tria del compasso} (1797): 
\begin{quote}
[\dots] mentre si trovano tante cose nuove progredendo nelle matema\-tiche, non si potrebbe forse trovare qualche luogo ancora incognito retrocedendo?\,\footnote{\translation{While we can find so many new things by moving forward in mathematics, why can't we find some still unknown place by retroceding?}}
\end{quote}
By respect for historical sources and for the reader's convenience, the quotations in the sequel will appear in the original tongue with a translation in square brackets, placed in footnote.

\section{The physical paradigm of mass-density \\versus the paradigm of primitives}\label{sez-paradigma}

In {\it  Philosophiae Naturalis Principia Mathematica} (1687) the first definition concerns mass and density:
\begin{quote}
 Quantitas materiae est mensura ejusdem orta ex ilius densitate et magni\-tudine conjunctim [\dots]. Hanc autem quantitem sub nomine corporis vel massa in sequentibus passim intelligo.\,\footnote{\translation{The quantity of matter is a measure of the matter itself, arising from its density and magnitude conjunctly [\dots]. It is this quantity that I mean  hereafter everywhere under the name of body or mass.}}
\end{quote}
In this sentence \textsc{Newton} presents the \emph{mass-density paradigm} (i.e.,  the mass can be computed in terms of the density and, conversely, the density can be obtained from the mass) as a fundament of Physics.

In \emph{Coexistent magnitudes} \cite[(1841)]{cauchy1841} \textsc{Cauchy}, with a clear didactic aim, uses the mass-density paradigm in order to give a unitary  exposition of several problems related to differential calculus.

From a mathematical point of view the implementation of this physical paradigm presents some difficulties and it does not assure a univocal answer. 
The first difficulty is in defining what is a ``mass'', 
the second is in choosing a procedure for evaluating  ``density'' and, 
finally, in determining under what condition and how it is possible ``to recover'' the mass from the density. 

All these critical aspects that we find in \textsc{Cauchy} \cite[(1841)]{cauchy1841}, are overcome in a precise and clear way by \textsc{Peano} in \emph{Applicazioni geometriche} \cite[(1887)]{peano87}. 

Natural properties that connect density and mass are the following:
\begin{enumerate}
\item \emph{The density of a homogenous body is constant.}\labelpag{hom}
\item \emph{The greater is the density, the   greater is the mass.} \labelpag{gre}
\item \emph{The mass of a body, as well as its volume, is the sum of its parts.}\labelpag{sum}
\end{enumerate}

The realization of  the physical paradigm can be mathematically expressed by the following formula
\begin{equation}\label{mass-volume}
\mu(A)=\int_A g \,\de({\rm vol}_n)
\end{equation}
where $\mu$ is the ``mass'', $g$ is the ``density'' and ${\rm vol}_n$ is the $n-$dimensional volume.\,\footnote{In today terminology, 
the realization of (\ref{mass-volume}) is expressed by saying that $g$ is the \emph{Radon-Nikodym derivative} of 
$\mu$ with respect to ${\rm vol}_n$.}

The properties \eqref{hom}, \eqref{gre} and \eqref{sum} do not allow for a direct derivation of \eqref{mass-volume} without further conditions depending on the meaning of integral;
for instance, having in mind the Riemann integral, an obvious necessary condition is the Riemann integrability of the density $g$.

In \textsc{Peano}'s \emph{Applicazioni Geometriche} \cite[(1887)]{peano87}:
\begin{itemize}
\item
the ``masses'' and the ``volumes'' are represented by \emph{distributive set functions}, as it will be shown in detail in \S \ref{sez-distributive}, 
\item
the ``densities'' (strict derivatives) are computed using a limit procedure, as we shall see in the sequel (see formula \eqref{der-f-peano}),
\item
the ``mass'' is recovered by integration using \eqref{mass-volume}. 
This final step is strengthened by the fact that \textsc{Peano}'s strict derivative is continuous.
\end{itemize}

The mathematical realization of mass-density paradigm is directly connected with mathematical paradigm of primitives, that is with the study of conditions assuring that integration is the inverse operation of differentiation.

At the beginning of the $20^{th}$ century the problem of looking for primitives is the cornerstone of the new theory of measure,  founded by \textsc{Lebesgue} \cite[(1904)]{lebesgue1904}. 
The problem of primitives  becomes arduous when one has to pass from functions of one variable to functions of more variables. 
\textsc{Lebesgue} in \emph{L'int\'egration des fonctions discontinues} \cite[(1910)]{lebesgue1910} overcomes these difficulties by substituting the integral of a generic function $g$   with  a set function $\mu$ described by formula \eqref{mass-volume}.

The paradigm of primitives gives more importance to the operations (of differentiation and integration) than to the set functions. On the contrary,  in the mass-density paradigm the primary aim is the evaluation of  the infinitesimal  ratio between two set functions
(for instance, mass and volume) in order to recover the ``mass'' by integrating the ``density'' with respect to ``volume''. On the other hand in the paradigm of primitives the main problem is an extension of the notion of integral in order to describe a primitive of a given function  and, consequently, to preserve fundamental theorem of calculus.

In \textsc{Lebesgue}'s works the two paradigms appear simultaneously for the first time in the second edition of his famous book {\em Le{\c c}ons sur l'int{\'e}gration et la recherche des fonctions
primitives} \cite[(1928) pp.\,196-198]{lebesgue1928}. In 1921 (see \cite[vol.\,I, p.\,177]{lebesgue_opere}) \textsc{Lebesgue} has already used some physical concept in order to make the notion of set function intuitive; analogously in  \cite[(1926)]{lebesgue1926} and  \cite[(1928) pp.\,290-296]{lebesgue1928} he uses the mass-density paradigm in order to make more natural the operations of differentiation and integration. In his lectures
\emph{Sur la mesure des grandeurs} \cite[(1935)]{lebesgue1935}, the physical paradigm leads \textsc{Lebesgue} to an alternative definition of derivative: he replaces his derivative of 1910 with the new uniform-derivative
(equivalent to the strict derivative introduced by \textsc{Peano}), thus allowing him to get continuity
of the derivative.

Before comparing \textsc{Peano}'s and \textsc{Lebesgue}'s derivative of set functions, we recall 
the definitions of derivatives given by \textsc{Peano} and \textsc{Cauchy}.

\textsc{Peano}'s strict derivative of a set function (for instance, the ``density'' of a ``mass'' $\mu$ with respect to the ``volume'') at a point $\bar x$ is computed, when it exists, as the limit of the quotient of the ``mass'' with respect to the ``volume'' of a cube $Q$, when the supremum of the distances of the points of the cube from $\bar x$ tends to 0 (in symbols $Q\to\bar x$). In formula, \textsc{Peano}'s strict derivative $g_{P}(\bar x)$ of a mass $\mu$ at $\bar x$ is given by:
\begin{equation}\label{der-f-peano}
 g_P(\bar x):=\lim_{Q\to \bar x}\frac{\mu (Q)}{\vol_n(Q)} \,\,\,.\end{equation}
Every limit procedure of a quotient of the form $\frac{\mu (Q)}{\vol_n(Q)}$ with $Q\to\bar x$ and the point $\bar x$ not necessarily belonging to $Q$, will be referred to as  \emph{derivative \`a la Peano}.

On the other hand, \textsc{Cauchy}'s derivative \cite[(1841)]{cauchy1841} is obtained as the limit between ``mass" and ``volume'' of a cube $Q$ {\it including} the point $\bar x$, when $Q\to\bar x$. In formula, \textsc{Cauchy}'s derivative $g_{C}(\bar x)$ of a mass $\mu$ at $\bar x$ is given by:
\begin{equation}\label{der-f-lebesgue} 
g_C(\bar x):=\lim_{\begin{subarray}{c}Q\to \bar x\\ \bar x\in Q\end{subarray}}\frac{\mu (Q)}{\vol_n(Q)} \,\,\,. \end{equation}
Every limit procedure of a quotient of the form $\frac{\mu (Q)}{\vol_n(Q)}$ with $Q\to\bar x$ and the point $\bar x$ belonging to $Q$, will be referred to as  \emph{derivative \`a la Cauchy}.

\textsc{Lebesgue}'s derivative of set functions is computed \emph{\`a la Cauchy}.
Notice that \textsc{Lebesgue} considers finite $\sigma$-additive and absolutely continuous measures as ``masses", while \textsc{Peano} considers distributive set functions. \textsc{Lebesgue}'s derivative exists (i.e., the limit $(2.6)$ there exists for \emph{almost every} $\bar x$), it is measurable and the reconstruction of a ``mass" as the integral of the derivative is assured by  absolute continuity
of the ``mass'' with respect to volume. On the contrary, \textsc{Peano}'s strict derivative does not
necessarily exist, but when it exists, it is continuous and the mass-density paradigm holds.\,\footnote{Clearly, if \textsc{Peano}'s strict derivative of a 
finite $\sigma$-additive measure exists, then  it coincides with \textsc{Lebesgue} derivative and the ``mass'' is absolutely continuous.

Nowadays it is not surprising that \textsc{Lebesgue}'s derivative can be seen as \textsc{Peano}'s strict derivative by lifting  measures on a $\sigma$-algebra $\Ac$ and $\Ac$-measurable functions to  measures on the Stone space  associated to $\Ac$ and the related continuous functions, respectively.}

The constructive approaches to differentiation of set functions  corresponding to the two limits \eqref{der-f-peano} and \eqref{der-f-lebesgue} are opposed to the approach given by \textsc{Radon} \cite[(1913)]{radon1913} and \textsc{Nikodym} \cite[(1930)]{nikodym1930}, who define the derivative in a more abstract and wider context  than those of \textsc{Lebesgue} and \textsc{Peano}. As in the case of \textsc{Lebesgue}, a Radon-Nikodym derivative exists; its existence is assured by assuming absolute continuity and 
$\sigma$-additivity of the measures.

In concluding this Section, let us remark that the physical properties
\eqref{hom}, \eqref{gre} and \eqref{sum}, that stand at the basis of
the mass-density paradigm, lead to the following direct characterization of the
Radon-Nikodym derivative. 
Let $\mu$ and $\nu$ be  finite $\sigma$-additive measures on a
$\sigma$-algebra $\Ac$ of subsets of $X$ and let $\nu$ be positive and
$\mu$ be absolutely continuous with respect to $\nu$.
A function $g:X\to\R$ 
is a \emph{Radon-Nikodym derivative} of $\mu$ with respect to $\nu$ (i.e., $\mu (A)=\int_A g\,\de\nu$ for every  $A\in\Ac$)
if and only if  the following two properties hold for every real number $a$:
\begin{enumerate}
\item $\mu(A)\geq a\,\nu(A)$ for every $A\subset \{g\geq a\}$ and $A\in\Ac$, \labelpag{Hahn_1}
\item $\mu(A)\leq a\,\nu(A)$ for every $A\subset \{g\leq a\}$ and $A\in\Ac$, \labelpag{Hahn_2}
\end{enumerate}
where $\{g\le a\}:=\{x\in X: g(x)\le a\}$ and, dually, $\{g\ge a\}:=\{x\in X: g(x)\ge a\}$. These properties \eqref{Hahn_1} and \eqref{Hahn_2}, expressed by \textsc{Nikodym} \cite[(1930)]{nikodym1930} in terms of Hahn decomposition of measures,  are a natural translation of properties \eqref{hom}, \eqref{gre} and \eqref{sum}.

\section{Peano on (pre-Lebesgue) classical measure theory}\label{sez-misura}

The interest of \textsc{Peano} in measure theory is rooted in his criticism of the \emph{definition of area} (1882), of the \emph{definition of integral} (1883) and of the \emph{definition of derivative} (1884).  This criticism  leads him to an innovative measure theory, which is extensively exposed in Chapter V of  \emph{Applicazioni geometriche} \cite[(1887)]{peano87}. 

The definition of area  given by \textsc{Serret} in \cite[(1879)]{serret} contrasted with the traditional definition of area: in 1882 \textsc{Peano}, independently of \textsc{Schwarz},  observed (see \cite[(1890)]{peano_area1890}) that the area of a cylindrical surface cannot be evaluated as the limit of inscribed polyhedral surfaces, as prescribed by \textsc{Serret}'s definition. In \emph{Applicazioni geometriche}, \textsc{Peano} provides a consistent definition of area and proves the integral formula of area.\,\footnote{This topic  will be extensively analyzed in a forthcoming paper by \textsc{Greco, Mazzucchi, Pagani}~\cite{gremazpagAREA}.}

\textsc{Peano}'s criticism of the definition of Riemann integral of a function and of its relation with the area of the \emph{ordinate-set} (i.e., hypograph of the function) \cite[(1883)]{peano1883}, forces him to introduce  outer/inner measure as the set-theoretic counterparts of upper/lower integral: he defines the latter in terms of infimum/supremum (instead of limits, as done traditionally) of the Darboux sums.\,\footnote{According with \textsc{Letta} \cite{letta}, the notion of negligible set  is introduced after an arduous process of investigation 
on ``similar'' notions related to cardinality and topology,
between 1870 and 1882. Afterward the definition of {\it Inhalt} (content) appears in the works by \textsc{Stolz} \cite[(1884)]{stolz1884}, \textsc{Cantor} \cite[(1884)]{cantor1884}, \textsc{Harnack} \cite[(1885)]{harnack1885}. 
The notions of inner and outer measure are introduced by \textsc{Peano} in  
\cite[(1883) p.\,446]{peano1883} and in
\cite[(1887)]{peano87}, and later by \textsc{Jordan} \cite[(1892)]{jordan1892}.
In the following we will refer to the inner and to the outer measures as to \emph{Peano-Jordan measures}.
}
\textsc{Peano}, in introducing the inner and outer measure as well as in defining area \cite[(1890)]{peano_area1890}, is also influenced by \textsc{Archimedes}'s approach on calculus of area, length and volume of convex figures.

In 1884, by analyzing the proof of mean  value theorem, given by \textsc{Jordan}
\footnote{\textsc{Jordan}, famous geometer and algebraist, publishes only a few papers on mathematical analysis. His most famous work is the {\it Cours d'analyse}, published in several editions. To our knowledge  the relationship between \textsc{Peano} and \textsc{Jordan} was good and based on reciprocal appreciation, as one can deduce from two letters conserved in Archives de la Biblioth\`eque Centrale de l'Ecole Polytechnique (Paris).}
in the first edition of \emph{Cours d'analyse}, \textsc{Peano} stresses the difference between differentiable functions and functions with continuous derivative. The continuity of derivative is expressed by \textsc{Peano} in terms of the existence of the limit
\begin{equation}\label{str-diff}\lim_{\begin{subarray}{c}x,y\to \bar x \\ x\neq y\end{subarray}}\frac{f(x)-f(y)}{x-y}\end{equation}
for any $\bar x$ in the domain of $f$.\,\footnote{
Later, in a paper with didactic value \cite[(1892)]{peano_der1892}, \textsc{Peano} re-proposes the distinction between  Definition \eqref{str-diff} and the usual derivative of a function, and underlines the correspondence of \eqref{str-diff} with the definition of density in Physics.

Nowadays the function $f$ is said \emph{strictly differentiable}  at the point $\bar x$  if the limit \eqref{str-diff} exists; consequently, the value of the limit $(3.1)$ is called \emph{strict derivative} of $f$ at $\bar x$.} 
Moreover, \textsc{Peano}, in his correspondence with \textsc{Jordan}
\cite[(1884)]{peano_jordan,peano_gilbert}, observes that uniform convergence of the difference quotient is equivalent to the continuity of the derivative.\,\footnote{
Section 80 of \textsc{Jordan}'s \emph{Cours d'analyse} \cite[(1893) p.\,68]{jordan1893}, titled ``\emph{Cas o\`u $\frac{f(x+h)-f(x)}{h}$ tend uniform\'ement vers $f'(x)$}'',
 contains a trace of it.} 
This notion of continuous derivative will be the basis of \textsc{Peano}'s strict derivative of distributive set functions.

\emph{Applicazioni geometriche}  is a detailed exposition (more than 300 pages) of several topics of geometric applications of infinitesimal calculus.\,\footnote{
As detailed in \textsc{Dolecki, Greco} \cite{greco-dolecki},
between several interesting concepts studied in \emph{Applicazioni geometriche} that are not directly connected with measure theory, we recall  the limit of sequences of sets (now called \emph{Kuratowski limits}), the introduction of the concept of differentiability of functions (nowadays called \emph{Fr\'echet differentiability}), the definition of tangent cone (nowadays called \emph{Bouligand cone}), the necessary condition of optimality (nowadays called \emph{Fermat conditions}) and a detailed study of problems of maximun and minimun.
 } 
In \emph{Applicazioni geometriche} \textsc{Peano} refounds the notion of Riemann integral by means of inner and outer measures\,\footnote{The simultaneous construction of inner and outer measure is the basis of the evolution of the theory leading to Lebesgue measure. Fortunately, \textsc{Carath\'eodory} \cite[(1914)]{caratheodory1914} and \textsc{Hausdorff} \cite[(1919)]{hausdorff1919} put an end to the intoxication due to the presence of inner measure, as \textsc{Carath\'eodory} writes: 
\begin{quote}
Borel and Lebesgue (as well as Peano and Jordan) assigned an outer measure $m^*(A)$ and an inner measure $m_* (A)$ to every point set $A$ [...]. The main advantage, however, is that the new definition 
[i.e., the exterior measure of Charath\'eodory] is independent of the concept of an \emph{inner measure} \cite[(2004) p.\,72]{edgar}.
\end{quote}}, and extends it to abstract measures.
 The development of the theory is based on solid topological and logical ground and on a deep knowledge of set theory. He introduces the notions of closure, interior and boundary of sets. 

\textsc{Peano} in \emph{Applicazioni geometriche} \cite[(1887)]{peano87}, and later \textsc{Jordan} 
in the paper \cite[(1892)]{jordan1892} 
and in the second edition of \emph{Cours d'Analyse} \cite[(1893)]{jordan1893}, develop the well known concepts of classical measure theory, namely, measurability, change of variables, fundamental theorems of calculus, with some methodological differences between them.\,\footnote{In a first paper of 
\textsc{Jordan} \cite[(1892)] {jordan1892} and in a more extensive way in his \emph{Cours d'analyse} \cite[(1893)]{jordan1893}, we find several \textsc{Peano}'s results. There are, however,  methodological differences between their approaches: \textsc{Peano} constructs his measure by starting from polygons, while \textsc{Jordan} considers (in the 2-dimensional case) squares. The definition proposed by \textsc{Peano} does not have the simplicity of that of \textsc{Jordan}, but it is independent of the reference frame and it is, by definition, invariant under isometries, without any need of further proof. Moreover, \textsc{Peano}'s definition allows for a direct computation of the proportionality factor appearing  under the action of affine transformation (in previous works \textsc{Peano} had developed a formalism allowing for computation of areas of polygons in a simple way, see \cite{gremazpagAREA}) for details. }

The mathematical tools employed by  \textsc{Peano} were really innovative at that time (and maybe are even nowadays), both on a geometrical and a topological level. \textsc{Peano} used extensively the geometric vector calculus introduced by \textsc{Grassmann}. The geometric notions include oriented areas and volumes (called \emph{geometric forms}). 

Our main interest concerns Chapter V of \textsc{Peano}'s \emph{Applicazioni geometriche}, where   we find differentiation of distributive set functions.

\emph{Applicazioni geometriche}  is widely cited, but we have the feeling that the work is not sufficiently known.
The  revolutionary character of \textsc{Peano}'s book is remarked by  J. \textsc{Tannery} \cite[(1887)]{tannery1887}:
\begin{quote}
Le Chapitre V porte ce titre: {\it Grandeurs g\'eom\'etriques}. C'est peut-\^etre le plus important et le plus int\'eressant, celui, du moins, par lequel le Livre de M. Peano se distingue davantage des Trait\'es classiques: les d\'efinitions qui se rapportent aux {\it champs de points}, aux points ext\'erieurs, 
int\'erieurs ou limites par rapport \a un champ, aux fonctions distributives (coexistantes d'apr\`es Cauchy), \a la longueur (\a l'aire ou au volume) externe, interne ou propre d'un champ, la notion d'int\'egrale \'etendue a un champ sont pr\'esent\'ees sous une forme abstraite, tr\`es pr\'ecise et tr\`es  claire.\,\footnote{\translation{%
Chapter V is titled: Geometric magnitudes. This chapter is probably the most relevant and
interesting, the one that marks the difference of the Book of Peano with respect to other
classical Treatises: definitions concerning \emph{sets of points}, exterior, interior and limit points
of a given set, distributive functions (coexistent magnitudes in the sense of Cauchy),
exterior, interior and proper length (or area or volume) of a set, the extension of the notion of integral
to a set, are stated in an abstract, very precise and very clear way.}}
\end{quote}

Only a few authors fully  realized the innovative value of Chapter V of
 \emph{Applicazioni geometriche}. As an instance, \textsc{Ascoli} says:
\begin{quote}
In [\emph{Applicazioni geometriche}] vi sono profusi, in forma cos\`i semplice da parere definitiva, idee e risultati divenuti poi classici, come quelli sulla misura degli insiemi, sulla rettificazione delle curve, sulla definizione dell'area di una superficie, sull'integrazione di campo, sulle funzioni additive di insieme; ed altri che sono tutt'ora poco noti o poco studiati [\dots].%
\,\footnote{
\translation{In \emph{Applicazioni geometriche} it is possible to find 
a clear and definitive exposition of many mathematical concepts and results, 
nowadays of common knowledge: results on measure of sets, on length of arcs, on the definition
of area of a surface, on the integration on a set, on additive set functions;
and other results that are not well known [\dots]
}  
}
\end{quote}

Most of the modern historians are aware of the  contributions to measure theory given by \textsc{Peano} and \textsc{Jordan} concerning inner and outer measure and measurability.\,\footnote{%In spite of some rare exceptions, this facts are well known.
To our knowledge the latest example of historian who forgot to quote any \textsc{Peano}'s contributions, is \textsc{Hochkirchen} \cite[(2003)]{hochkirchen2003}. Ironically, the symbols $\underline\int$ and $\overline\int$ which \textsc{Volterra} introduced for denoting lower and upper integral, were ascribed to \textsc{Peano} by \textsc{Hochkirchen}.
} 

Only a few historians mention \textsc{Peano}'s contributions to derivative of set functions: \textsc{Pesin} \cite{pesin}, \textsc{Medvedev} \cite{medvedev} and \textsc{Hawkins}  \cite{hawkins} and others. 

\textsc{Pesin} \cite[(1970)\, pp.\,32-33]{pesin}, who does ``not intend to overestimate the importance of \textsc{Peano}'s results'',   recalls some results of \textsc{Peano}'s work without giving details or appropriate definitions.

\textsc{Medvedev}  in \cite[(1983)]{medvedev} recalls \textsc{Peano}'s contributions giving detailed information  both on the integral as a set function and on the \textsc{Peano}'s derivative. In our opinion he gives an excessive importance to mathematical priorities without pointing out the differences between \textsc{Peano}'s contribution of 1887  and \textsc{Lebesgue}'s contribution of 1910.\,\footnote{
\textsc{Dieudonn\'e}, reviewing in \cite[(1983)]{dieudonne} the \textsc{Medvedev}'s paper \cite[(1983)]{medvedev}, with his usual sarcasm denies any logical value of \textsc{Peano}'s definitions concerning limits and sets.
Against any historical evidence, \textsc{Dieudonn\'e} forgets several \textsc{Peano}'s  papers
on several notions of limit, and ignores the \emph{Formulario mathematico} where \textsc{Peano}
presents a large amount of mathematical results, including set axiomatization, through his logical \emph{ideography}. Besides, \textsc{Dieudonn\'e} forgets \textsc{Bourbaki}'s \emph{Elements of the history of mathematics} and ignores that  the building blocks of \textsc{Peano}'s ideography  are  the atomic propositions: $x\in X$ and $x=y$.}

\textsc{Hawkins} does not describe  \textsc{Peano}'s results on differentiation and integration in detail, as they are too far from the main aim of  his book, but he is aware of   \textsc{Peano}'s contributions to differentiation of set functions \cite[p.\,88,\,185]{hawkins},  and appraises \textsc{Peano}'s book  \emph{Applicazioni geometriche}:
\begin{quote}
the theory is surprisingly elegant and abstract for a work of 1887 and strikingly modern in his approach \cite[p.\,88]{hawkins}.
\end{quote} 

None of the historian quoted above, establishes a link between \textsc{Peano}'s work on differentiation of measure in  \emph{Applicazioni geometriche} with his paper \emph{Grandezze coesistenti} \cite{peano1915} and with Lebesgue's comments on differentiation presented in \emph{La mesures des grandeurs} \cite[(1935)]{lebesgue1935}.

Main primary sources on which our paper is based are \cite{cauchy1841,peano1915,lebesgue1910,lebesgue1935,vitali1915,vitali1916,fubini1915b,fubini1915a}. 

\section{Cauchy's coexistent magnitudes}\label{sez-cauchy}
\textsc{Cauchy}'s seminal paper \emph{Coexistent magnitudes} \cite[(1841)]{cauchy1841} presents some difficulties for the modern reader: the terms he introduces are rather obscure (for instance, \emph{grandeurs}, \emph{coexistantes}, \emph{\'el\'ements}, \dots), and the reasonings are based on vague geometric language, accordingly to the \textsc{Cauchy}'s taste. Actually, \textsc{Cauchy}'s aim was to make mathematical analysis as well rigorous as geometry \cite[(1821) p.\,ii]{cauchy1821}:
\begin{quote}
Quant aux m\'ethodes, j'ai cherch\'e \a leur donner toute la rigueur qu'on exige en g\'eom\'etrie, de mani\`ere \`a ne jamais recourir aux raisons tir\'ees de la g\'en\'eralit\'e de l'alg\`ebre.\,\footnote{\translation{%
About methods, I have tried to be rigorous as required in geometry, 
in order to avoid the general reasonings occurring in algebra.}

Not all mathematicians at that time considered geometry as a model of rigor. Indeed  
\textsc{Lobachevsky} starts his famous book  ``Theory of parallels'' \cite[(1829) p.\,11]{Loba1829} with the following sentence:
\begin{quote}
In geometry I find certain imperfections which I hold to be the reason why this science, apart from transition into analytics, can as yet make no advance from that state in which it has come to us from Euclid.
\end{quote}}
\end{quote}

In his \emph{Le\c cons de m\'ecanique analytique} \cite[(1868) pp.\,172-205]{moigno1868}
\textsc{Moigno}, a follower of \textsc{Cauchy}, reprints the paper \emph{Coexistent magnitudes}. He puts into evidence the
vagueness of some terms of \textsc{Cauchy}, unfortunately without adding any comment that may help the reader to a better understanding of \textsc{Cauchy}'s paper itself. 

The meaning of the terms  ``\emph{grandeurs}'' and ``\emph{coexistantes}'' can be made precise by analyzing the list of examples given by \textsc{Cauchy}.
 He implicitly postulates the following properties of  ``\emph{grandeurs}'':
\begin{enumerate}
\item \labelpag{cond2} a magnitude can be divided into finitely many infinitesimal equal elements
(using the terminology of \textsc{Cauchy}), where infinitesimal is related to magnitude 
and diameter; 
\item  the ratio between coexistent magnitudes (not necessarily homogeneous) is a numerical quantity.
\end{enumerate}
Concerning the term ``\emph{coexistantes}'', coexistent magnitudes are defined by \textsc{Cauchy} as ``magnitudes which exist together, change simultaneously and the parts of one magnitude exist and change 
in the same way as the parts of the other 
  magnitude''.\,\footnote{
\textsc{Cauchy} says in \cite[(1841) p.\,188]{cauchy1841}:
\begin{quote}
Nous appellons \emph{grandeurs} ou quantit\'es \emph{coexistantes} deux grandeurs ou
quantit\'es qui existent ensemble et varient simultan\'ement, de telle sorte
que les \'el\'ements de l'une existent et varient, ou s'\'evanouissent, en m\^eme temps que
les \'el\'ements de l'autre.
\end{quote}%
} 
Despite of the vagueness of this definition,
the meaning of ``\emph{coexistantes}" is partially clarified by many examples of coexistent magnitudes given by \textsc{Cauchy}
\cite[(1841) pp.\,188--189]{cauchy1841}, such as the volume and the mass of a body, the time and 
the displacement of a moving point, the radius and the surface of a circle, the radius and the volume of a sphere, the height and the area of a triangle, the height and the volume of a prism, the base and the volume of a cylinder, and so on.

Vagueness of the \textsc{Cauchy}'s definition of ``\emph{grandeurs coexistantes}'' was pointed out by \textsc{Peano}. 
In \emph{Applicazioni geometriche} \cite[(1887)]{peano87} and in \emph{Grandezze coesistenti} \cite[(1915)]{peano1915}, \textsc{Peano} defines them as set functions 
over the same given domain, satisfying  additivity properties in a suitable sense.

The primary aim of \textsc{Cauchy} is pedagogic: he wants  to write a paper making easier the study of infinitesimal calculus and its applications. As it is easy to understand, \textsc{Cauchy} bases himself on the mass-density paradigm and introduces  the limit of the average of two coexistent magnitudes, calling it  {\it differential ratio}. In a modern language we could say that the coexistent magnitudes are set functions, while the differential ratio is a point  function. \textsc{Cauchy} points out that the differential ratio is termed 
in different ways depending on the context, namely, on the nature of the magnitudes themselves (for instance, mass density of a body at a given point,  velocity of a moving point at a given time, hydrostatic pressure at a point of a given surface, \dots). 

Now we list the most significant theorems that are present in the paper of \textsc{Cauchy},
preserving, as much as possible, his terminology. 

\begin{theorem}\label{mediaintegrale}\cite[Theorem 1, p.\,190]{cauchy1841} The average between two coexistent magnitudes is bounded between the supremum  and the infimum of the values of the differential ratio.
\end{theorem}

\begin{theorem}\label{densitanulla}\cite[Theorem 4, p.\,192]{cauchy1841}  A magnitude vanishes whenever its differential ratio, with respect to another coexistent magnitude, is a null function.
\end{theorem}
\begin{theorem}\label{teomedia} \cite[Theorem 5, p.\,198]{cauchy1841} 
If the differential ratio between two coexistent magnitudes is a continuous function, then the ``mean value property'' holds.\,\footnote{Let $\mu,\nu : \Ac \to \R$ be two magnitudes and let $g$ be the differential ratio of $\mu$ with respect to $\nu$. We say that the \emph{mean value property} holds if, for any set $A \in \Ac$,
with $\nu(A) \ne 0$, there exists a point $P\in A$ such that $g(P)=\frac{\mu(A)}{\nu (A)}$.}
\end{theorem}
\begin{theorem}\label{densitauguale}  \cite[Theorem 13, p.\,202]{cauchy1841} If two magnitudes have the same differential ratio with respect to another magnitude, then they are equal.
\end{theorem}

Even if \textsc{Cauchy} presents proofs that are rather ``vanishing'', his statements (see theorems listed above)
and his use of the differential ratio allow \textsc{Peano} to rebuild
his arguments on solid grounds. 
\textsc{Peano} translates the coexistent magnitudes into the concept of distributive set functions, restating the theorems presented by \textsc{Cauchy} and proving them rigorously.

In \textsc{Peano}, the property of continuity of the differential ratio (whenever it exists) is a 
consequence of its definition. 
On the contrary, \textsc{Cauchy}'s definition of differential ratio does not guarantee its continuity.
\textsc{Cauchy} is aware of the fact that the differential ratio can be discontinuous, nevertheless he thinks that, in the most common ``real'' cases, it may be assumed to be continuous; see \cite[(1841), p.\,196]{cauchy1841}: 
\begin{quote} 
Le plus souvent, ce rapport diff\'erentiel sera une fonction continue de la variable dont il d\'epend, c'est-\`a-dire qu'il changera de valeur avec elle par degr\'es insensibles.\,\footnote{\translation{%
Almost always, this differential ratio is a continuous function of the independent variable,
i.e., its values change in a smooth way.}}
\end{quote} 
and \cite[(1841) p.\,197]{cauchy1841}:
\begin{quote}
Dans un grand nombre de cas, le rapport diff\'erentiel $\rho$ est une fonction continue [\dots].\,\footnote{\translation{%
Almost always, the differential ratio $\rho$ is a continuous function [\dots].}}
\end{quote}

In evaluating the differential ratio as a ``limit of average values $\frac{\mu(A)}{\nu(A)}$ at a point $P$'', 
for \textsc{Peano} the set $A$ does not necessarily include the point $P$, 
while for \textsc{Cauchy} $A$ includes $P$ (as \textsc{Cauchy} says: $A$ \emph{renferme le point} $P$).

This difference is fundamental also in case of linearly distributed masses. Indeed a linear mass distribution, described in terms of a function of a real variable, 
admits a differential ratio in the sense of \textsc{Peano} if the derivative exists and is continuous,
whilst it admits a
differential ratio in the sense of \textsc{Cauchy}
\footnote{
Using the identity
\begin{equation*}
\frac{f(x+h) - f(x-k)}{h+k} = 
\frac{h\frac{f(x+h)-f(x)}{h} +k\frac{f(x)-f(x-k)}{k}}{h+k} 
\quad \text{ for every } k, h > 0
\end{equation*}
the reader can easily verify that the differential ratio in the sense of \textsc{Cauchy} 
exists (i.e., the limit of $\frac{f(x+h) - f(x-k)}{h+k}$ exists for $k\to 0^+$ and
$h \to 0^+$, with $h+k > 0$) whenever $f'(x)$ exists.
}
only if the function is differentiable 
\cite[(1841) p.\,208]{cauchy1841}:
\begin{quote}
Lorsque deux grandeurs ou quantit\'es coexistantes se r\'eduisent \`a une 
variable $x$ et \`a une fonction $y$ de cette variable, le rapport diff\'erentiel de
fonction \`a la variable est pr\'ecis\'ement ce qu'on nomme la \emph{d\'eriv\'ee} de la
fonction ou le \emph{coefficient diff\'erentiel}.\,\footnote{\translation{%
When two coexistent magnitudes are a variable $x$ and a function $y$ of $x$,
the differential ratio of the function with respect to the variable $x$ 
coincides with the \emph{derivative} of the function.}}
\end{quote}

Concerning the existence of the differential ratio, \textsc{Cauchy} is rather obscure; indeed  whenever he defines the differential ratio, he specifies that ``it will converge in general to a certain limit different from $0$''.  As \textsc{Cauchy} does not clarify the meaning of the expression ``in general'', 
the conditions assuring the existence of the differential ratio are not given explicitly.  
On the other hand, \textsc{Cauchy} himself is aware of this lack, as in several theorems he explicitly assumes that the differential ratio is ``completely determined at every point''.

Concerning the mass-density paradigm, in \textsc{Cauchy}'s \emph{Coexistent magnitudes} an explicit  formula allowing for constructing the mass 
of a body in terms of its density is also lacking. 
In spite of this, \textsc{Cauchy} provides a large amount of theorems and corollaries giving an approximate calculation of the mass under the assumption of continuity of the density.
We can envisage this approach as a first step toward the modern notion of integral with respect to 
a general abstract measure.

We can summarize further \textsc{Cauchy}'s results into the following  theorem:
\begin{theorem}\cite[(1841) pp.\,208--215]{cauchy1841} Let us assume that the differential ratio $g$ between two coexistent magnitudes $\mu$ and $\nu$ exists and is continuous. Then $\mu$ can be computed in terms of the integral of $g$ with respect to $\nu$.
\end{theorem}

\textsc{Cauchy} concludes his memoir \cite[(1841) pp.\,215--229]{cauchy1841} with a second section
in which he states the following theorem in order to evaluate lengths, areas and volumes of homothetic elementary figures.

\begin{theorem}\label{teo-prop} \cite[Theorem 1, p.\,216]{cauchy1841} Two coexistent magnitudes are proportional, whenever to equal parts of one magnitude there correspond equal parts of the 
other.\,\footnote{One can observe that this theorem holds true by imposing condition \eqref{cond2}.
}
\end{theorem}

Even if the \textsc{Cauchy}'s paper contains several innovative procedures, to our knowledge only a few authors (\textsc{Moigno}, \textsc{Peano}, \textsc{Vitali}, \textsc{Picone} and \textsc{Lebesgue}) quote it, and only \textsc{Peano} and \textsc{Lebesgue}  analyze it in details. 

\section{Distributive families, decompositions and Peano additivity}\label{sez-distributive}

In his paper \emph{Le grandezze coesistenti} \cite[(1915)]{peano1915},
\textsc{Peano} introduces a general concept of  \emph{distributive function}, namely 
a function $f:A\to B$, where $(A,+), (B,+)$ are two sets endowed with binary operations, denoted by the same symbol $+$, satisfying the equality
\begin{equation}
\label{distr}f(x+y)=f(x)+f(y) 
\end{equation}
for all $x,y$ belonging to $A$ and, if necessary, verifying suitable assumptions.\,\footnote{Among distributive functions considered by \textsc{Peano}, there are the usual linear functions 
and particular set functions. The  reader has to pay attention in order to avoid the interpretation of distributive set functions as finitely additive set functions.}
\textsc{Peano} presents several examples of distributive functions. 
As a special instance, $A$ stands for the family $\Po(X)$
of all subsets of a finite dimensional Euclidean space $X$,
``$+$'' in the left hand side of (\ref{distr}) is the union operation, 
and ``$+$'' in the right hand side of (\ref{distr}) is the logical OR (denoted in \textsc{Peano}'s
ideography by the same symbol of set-union); therefore, equation \eqref{distr} becomes: 
\begin{equation}\label{distr2}f(x\cup y)=f(x) \cup f(y).\end{equation}
 
To make (\ref{distr2}) significant, \textsc{Peano} chooses a family ${\mathcal U}\subset \Po(X)$
and defines ``$f(x)$'' as ``$x\in \mathcal{U}$''.
Consequently (\ref{distr2}) becomes: 
\begin{equation}\label{griglia}x\cup y\in {\mathcal U}\Longleftrightarrow x\in {\mathcal U}\text{ or  } y\in {\mathcal U} \end{equation}
for all $x,y \in \Po(X)$.
A family $\mathcal U$ satisfying (\ref{griglia}) 
is called by \textsc{Peano} a \emph{distributive family}.\,\footnote{This notion of distributive family will be  rediscovered later by \textsc{Choquet} \cite[(1947)]{choquet}, who called it   {\it grill} and recognized it as the dual notion of \textsc{Cartan}'s {\it filter} \cite[(1937)]{cartan}.}  

Moreover, \textsc{Peano} considers \emph{semi-distributive} families $\Fo\subset \Po (X)$, i.e., families of sets such that 
\begin{equation}
x\cup y\in \Fo \Longrightarrow  x\in \Fo \text{ or  } y\in \Fo
\end{equation} 
for all $x, y \in \Po(X)$.

A distributive family of subsets of $X$ is obtained by a semi-distributive family $\Fo$ by adding to $\Fo$ any supersets of its elements.
\textsc{Peano} states the following theorem, and attributes to \textsc{Cantor} \cite[(1884) p.\,454]{cantor1884} both its statement
and its proof.
\begin{theorem}[Cantor compactness property]\label{teo-F1} 
Let $\Fo$ be a semi-distributive family of a finite-dimensional Euclidean space, and
let $S$ be a bounded non-empty set belonging to $\Fo$. Then there exists a point $\bar x$, belonging to the closure of $S$, such that any neighborhood of $\bar x$ contains a set belonging to $\Fo$. 
\end{theorem}

The notion of distributive family is essential in the study of the derivation of distributive set functions by \textsc{Peano}.
Distributive families have been introduced by \textsc{Peano} in \emph{Applicazioni geometriche} in 1887. Moreover, he uses them in his famous paper on the existence of solutions of differential equations \cite[(1890) pp.\,201--202]{peano1890} and, later, in his textbook \emph{Lezioni di analisi infinitesimale} \cite[(1893) vol.\,2, pp.\,46--53]{peano1893}. The role played by this notion is nowadays recovered by ``compactness by coverings'' or by ``existence of accumulation points''.\,\footnote{Two examples of distributive families considered by \textsc{Peano} are
${\mathcal U}:=\{ A\subset \R^n\, :\; \text{ card}( A)=\infty  \}$,
 and
${\mathcal U}_h:=\{A\subset \R^n\, :\, \sup_A h=\sup_{\R^n}h\}$,
where $h:\R^n\to\R$ is a given real function.
} 

In proving Theorem \ref{teo-F1}, \textsc{Peano} decomposes a subset of the Euclidean space $\R^n$ following a grid of n-intervals 
implemented by cutting  sets along hyperplanes parallel to coordinate axis.
We may formalize this procedure in the following way.

Let us denote by $H$ a hyperplane of the form $H:=\{x\in\R^n\,:\, \langle x,e_i\rangle =a\}$ where $e_i$ is a vector of the canonical basis of $\R^n$ and $a\in\R$. Let us denote by $H^+$ and $H^-$ the two closed half-spaces delimited by $H$.

A family  $\Fo$ of subsets of $\R^n$ is called \emph{semi-distributive by cutting along hyperplanes} 
if 
$$   A\cap H^+ \in \Fo \text{ or } A\cap H^-\in \Fo $$
for every $A \in \Fo$ and for every hyperplane $H$ of $\R^n$ of the form indicated above.
Under this restrictions  a new version of Theorem \ref{teo-F1} still holds:

\begin{theorem}[Cantor compactness property by interval-decompositions]\label{teo-F2} 
Let $\Fo$ be semi-distributive by cutting along hyperplanes and
let $S$ be a bounded non-empty set belonging to $\Fo$. Then there exists a point $\bar x$ belonging to the closure of $S$ such that any neighborhood of $\bar x$ contains a set belonging to $\Fo$. 
\end{theorem}

To express additivity properties of set functions, \textsc{Peano}, as it was common at his time \footnote{A similar expression is used also by \textsc{Jordan} \cite{jordan1892}:
\begin{quotation}
[C]haque champ $E$ a une \'etendue d\'etermin\'ee; [\dots]  si on le d\'ecompose en plusieurs parties $E_1$, $E_2$, \dots, la somme des \'etendues de ces parties est \'egale \a l'\'etendue totale de $E$.
\translation{%
Every set $E$ has a defined extension; [\dots] if $E$ is decomposed into parts
$E_1$, $E_2$, \dots, the sum of the extensions of these parts is equal to the
extension of $E$.}
\end{quotation}

}, uses the term  \emph{decomposition}. \textsc{Peano} writes in \emph{Applicazioni geometriche} 
\cite[(1887) p.\,164, 167]{peano87}:
\begin{quote}
Se un campo $A$ \`e decomposto in parti $A_1,A_2, \dots, A_n$ esso si dir\`a \emph{somma} delle sue parti, e si scriver\`a 
$$A=A_1+A_2+\dots+A_n$$
[\dots] Una grandezza dicesi \emph{funzione distributiva} di un campo, se il valore di quella grandezza corrispondente ad un campo \`e la somma dei valori di essa corrispondenti alle parti in cui si pu\`o decomporre un campo dato.\,\footnote{\translation{%
If a set $A$ is decomposed into the parts $A_1,A_2,\dots,A_n$, it will be called \emph{sum}
of its parts, and it will be denoted by $A=A_1+A_2+\dots+A_n$. [\dots] A magnitude is said to be
a \emph{distributive set function} if its value on a given set is the sum of the corresponding
values of the function on the parts decomposing the set itself.}}
\end{quote}

In order to formalize in modern language both the operation of ``decomposing'' and his use in \textsc{Peano}'s
works, we can pursuit 
a ``minimal'' way, leading to ``families of interval-decompositions'', and a ``proof-driven'' way, 
leading to ``families of finite decompositions''. 

First, the minimal way consists in implementing the procedure of decomposing by cutting along hyperplanes used by \textsc{Peano} in proving Theorem \ref{teo-F1}.
More precisely, let $\Ac$ be a family of subsets of the Euclidean space $\R^n$; a finite family 
$\{A_i\}_{i=1}^m$  of elements of $\Ac$ is called
an \emph{interval-decomposition} of $A \in \Ac$ if it is obtained by iterating the procedure
of cutting by hyperplanes.
In other words, an interval-decomposition $\{A_i\}_{i=1}^m$ of a set $A$
is a finite sub-family of $\Ac$ defined recursively as follows:
\begin{itemize}
\item
for $m=1$, $A_1=A$; 
\item
for $m=2$, there exists a hyperplane $H$ such that $A_1= A\cap H^-$ and $A_2= A\cap H^+$;
\item
for $m>2$, there exist two distinct indices $i_0, i_1 \le n$ such that 
$\tilde A := A_{i_0} \cup A_{i_1} \in \Ac$ and the families
 $\{A_i : 1 \le i \le m, i \ne i_0, i \ne i_1 \}\cup \{\tilde A \}$ and $\{A_{i_0}, A_{i_1} \}$ are  interval-decompositions of $A$ and  $\tilde A$, respectively.
\end{itemize}
The totality of these interval-decompositions will be denoted by $\Dbb_{\textrm{int}}(\Ac)$. In the case where $\Ac$ is the family of all the closed bounded subintervals of a given closed interval $[a,b]$
of the real line,  
an arbitrary interval-decomposition of an interval $[a',b']\subset[a,b]$ is a family $\{ [a_{i-1},a_i] \}_{i=1}^m$
where  $a' = a_0 \le a_1 \le \dots \le a_{m-1} \le a_m= b'$. 
 The totality of these interval-decompositions are denoted by
$\Dbb_{\textrm{int}}(a,b)$.

The second way  
consists in summarizing explicitly the properties of the decompositions themselves, as used by \textsc{Peano}
in defining the integral and in proving related theorems
\footnote{
See pages 165 and 186-188 of \emph{Applicazioni geometriche} \cite[(1887)]{peano87}. 
}, as it will be seen in Section \ref{sez-massdensity}. This leads to the following definitions of \emph{family of finite decompositions} and of the related \emph{semi-distributive family}, \emph{Cantor compactness property} and
\emph{distributive set functions}.

Let $\Ac$ be again a family of subsets of an Euclidean space $\R^n$ and let us denote by $\Po_{f}(\Ac)$ the set  of all non-empty finite subfamily of $\Ac$.
Define $\Ubb (\Ac)$ by
$$\Ubb (\Ac):=\{\Hc\in \Po_{f}(\Ac): \cup \Hc \in \Ac\}.
$$
Let $\Dbb$ be a subset of $\Ubb(\Ac)$;
we will say that $\Hc$ is a $\Dbb$-\emph{decomposi\-tion} of $A$ if $\Hc \in \Dbb$ and  $A=\cup\Hc$.

\begin{definition}
$\Dbb  \subset \Ubb(\Ac)$ is called a \emph{family of finite decompositions 
relative to $\Ac$} if  the following properties are satisfied: 
\begin{enumerate}
\item $\{A\}\in \Dbb$ for every $A\in\Ac$;
\item \labelpag{dec_1} if $\Hc$ and $\Gc$ are $\Dbb$-decompositions of a set $A$, then 
$$\{H\cap G: H\in\Hc,\,G\in\Gc 
\}$$ is  a 
$\Dbb$-decomposition of $A$; 
\item \labelpag{dec_2} if $\Hc$ and $\Gc$ are $\Dbb$-decompositions of $A$, then for every $G\in \Gc$
the family 
$$\Hc_G := \{H\cap G : H\in\Hc 
\}$$ 
 is a $\Dbb$-decomposition of $G$;
\item \labelpag{dec_3} if $\Hc$ is a $\Dbb$-decomposition of $A$ and, moreover, for every  $H \in \Hc$ the family 
$\Gc_H$ is a $\Dbb$-decomposition of $H$, then 
$$\cup \{\Gc_H : H \in \Hc\}$$
 is a  
$\Dbb$-decomposition of $A$.
\end{enumerate}
\end{definition}

\begin{definition} \label{def_inf}
A family $\Dbb$ of finite decompositions relative to $\Ac$ is called \emph{infinitesimal} if,
for every  bounded set $A \in \Ac$ and for every real number $\varepsilon >0$, there is
a $\Dbb$-decomposition $\Hc$ of $A$ such that the diameter of every $H \in \Hc$
is less than $\varepsilon$. 
\end{definition}

\begin{definition}
Let $\Dbb$ be a family of finite decompositions 
relative to $\Ac$. Then 
a \emph{set function} $\mu \colon \Ac \to \R$ is said to be \emph{distributive with respect to} $\Dbb$, if 
\begin{equation}\mu (\cup \Hc) = \sum_{H \in \Hc} \mu(H) \text{ for every }\Hc \in \Dbb.\end{equation}
\end{definition}
Consequently,  
\begin{definition} \label{def_semi}
Let $\Dbb$ be a family of finite decompositions 
relative to $\Ac$. A family $\Fo$ of subsets of the Euclidean space $\R^{n}$ is  said to be \emph{semi-distributive with respect to} $\Dbb$, if 
\begin{equation}\Hc\in\Dbb \text{ and } \cup\Hc\in \Fo \Longrightarrow \exists H\in\Hc \text{ such that } H\in\Fo.
\end{equation}

\end{definition}

\begin{theorem}[Cantor compactness property by an arbitrary family of decompositions]\label{teo-F3} 
Let $\Dbb$ be an infinitesimal family of finite decompositions 
relative to $\Ac$ and  let $\Fo$ be a semi-distributive family  with respect to $\Dbb$.  If  $S$ is a bounded non-empty set belonging to $\Fo$. Then there exists a point $\bar x$ belonging to the closure of $S$ such that any neighborhood of $\bar x$ contains a set belonging to $\Fo$. 
\end{theorem}

In the following, an expression of type ``$\mu \colon(\Ac,\Dbb) \to \R$ is a distributive set function''
stands for ``$\Dbb$ is a family of finite
decompositions relative to $\Ac$ and $\mu\colon\Ac \to \R$ is a distributive set function with respect to  $\Dbb$. 

Examples of families of decompositions are $\Ubb(\Ac)$, and 
\begin{enumerate}
\item \labelpag{dec-inter}
the family $\Dbb_{\textrm{int}}(\Ac)$ of all  interval-decompositions introduced above;\,\footnote{Notice that
a real valued set function $\mu:\Ac \to \R$ is distributive with respect to $\Dbb_{\textrm{int}}(\Ac)$, if $\mu (A)=\mu (A\cap H^+)+\mu (A\cap H^-)$
for every $A\in \Ac$ and for every hyperplanes $H$ such that
$A\cap H^+ \in \Ac $ and $A\cap H^-\in \Ac$. Inner and upper Peano-Jordan measures are both distributive in this sense, but they are not finitely additive. 
}
\item \labelpag{dec-picone}
the family of all $\Hc \in \Ubb(\Ac)$ 
such that the interiors of two arbitrary distinct elements of $\Hc$ have empty intersection
and every $H \in \Hc$ is Peano-Jordan measurable;
\item \labelpag{dec-picone2}
the family of all $\Hc \in \Ubb(\Ac)$ 
such that the intersection of the closure of two arbitrary distinct elements of $\Hc$ have  null  Peano-Jordan measure and every  $H\in \Hc$ is bounded;
\item \labelpag{last_ex}
the family of all $\Hc \in \Ubb(\Ac)$ 
such that two arbitrary distinct elements of $\Hc$ have empty intersection.

\end{enumerate}  

The interval-decompositions (in particular $\Dbb(a,b)$) occurs frequently in
\textsc{Peano}'s works. 
Distributive set functions related to the last example (\ref{last_ex})  are well
known as \emph{finitely additive set functions}; this type of additivity, expressed in terms of partitions of sets, was introduced for the first time in \textsc{Borel} \cite[(1898), pp.\,46-50]{borel1898}, and, more clearly, in \textsc{Lebesgue} \cite[(1902), p.\,6]{lebesgue1902}. 

As far as we know, all historians interpreted \textsc{Peano}'s distributive set functions
 as ``finitely additive'' set functions.\,\footnote{Observe that  inner and  outer Peano-Jordan measures on  Euclidean spaces  are not finitely additive, but they  are  distributive set functions with respect to the families of decomposition of type (\ref{dec-inter}) or (\ref{dec-picone}). Moreover, notice that  outer Peano-Jordan measure is  a distributive set function with respect to a family of decompositions of type  (\ref{dec-picone2}).} 
For instance, in the proof of the integrability of functions \cite[(1887) p.\,188]{peano87}, \textsc{Peano} uses  distributivity properties of the upper and lower integral with respect to the domain of integration; clearly neither
the upper nor the lower integral are finitely additive.

\section{Peano's strict derivative of distributive functions \\ and its applications}\label{sez-derivata}

In \emph{Applicazioni geometriche} \cite[(1887)]{peano87} \textsc{Peano} translates in terms of ``distributive functions'' the ``magnitudes'' of \textsc{Cauchy}, so that two \textsc{Cauchy}'s magnitudes are ``coexistent'' if they 
are distributive functions with the same domain.

\textsc{Peano}'s distributive set functions are called \emph{positive} if their values are positive.
\textsc{Peano}'s \emph{strict derivative} is defined by
\footnote{
In \textsc{Peano}'s words \cite[(1987) p.\,169]{peano87}:
\begin{quote}
Diremo che, in un punto $P$, il \emph{rapporto} delle due funzioni distributive $y$ ed $x$
d'un campo vale $\rho$, se $\rho$ \`e il limite verso cui tende il rapporto dei valori
di queste funzioni, corrispondenti ad un campo di cui tutti i punti si avvicinano
indefinitamente a $P$. 
\end{quote}
\translation{%
Given two distributive functions 
$y$ an $x$ defined over a given set, we say that their \emph{ratio}, at a given point $P$, is $\rho$, 
if $\rho$ is the limit of the ratio between the values of the two functions, taken along
sets for which all its points approach the point $P$.}}
\begin{definition}
Let $\mu, \nu : (\Ac,\Dbb) \to \R$ be distributive set functions, 
and let $\nu$ be positive. 
A real function $g$ over a set $S$ is called a ``strict derivative of $\mu$ with respect to $\nu$'' on $S$
(denoted by $\frac{d\,\mu}{d\,\nu}$ and termed \emph{rapporto} in \emph{Applicazioni geometriche})
if, for every point $x \in S$ and
for every $\epsilon >0$, there exists  $\delta >0$ such that 
\footnote{One can note that for the definition of strict derivative at a point $x$, the point $x$ itself
must be an accumulation point with respect to the family $\Ac$ and the measure $\nu$, that is,
for all $\delta >0$, there exists a $A\in\Ac$ such that $\nu(A)\neq 0$ and $A\subset B_\delta (x)$, where
$B_\delta (x)$ denotes the Euclidean ball of center $x$ and radius $\delta$.}
\begin{equation}
\left|\frac{\mu (A)}{\nu (A)}-g(x)\right|<\epsilon \quad \text{for every } A\in \Ac, \, \text{with } \nu (A)\neq 0, \, 
A\subset B_\delta(x).
\end{equation}
\end{definition}

It is worth noticing that the concept of strict derivative given by \textsc{Peano} provides a consistent mathematical
ground to the concept of ``infinitesimal ratio'' between two magnitudes,
successfully used since \textsc{Kepler}.
A remarkable example given by \textsc{Peano} is the evaluation of a rectifiable  arc length by integrating the ``infinitesimal
arc length'' $ds$. Notice that, whenever $ds$ exists in the sense of \textsc{Peano}, the corresponding
integral provides the length of the arc. On the contrary, the integration of
the infinitesimal arc length $ds$, evaluated in the sense of \textsc{Lebesgue} (1910), provides  the length of the arc
only in case of absolute continuity of the arc parametrization (see \textsc{Tonelli} \cite[(1908)]{tonelli}) 

The existence of \textsc{Peano}'s strict derivative is not assured in  general; 
its characterizing properties are clearly presented in \emph{Applicazioni geometriche} and can be summarized in the following theorems.

First, \textsc{Peano} gives a precise form to \textsc{Cauchy}'s Theorem \ref{mediaintegrale}, stating the following:

\begin{theorem}[see \textsc{Peano} {\cite[Theorem 13, p.\,170]{peano87}}  for $\Dbb=\Dbb_{\textrm{int}}$]\label{teorema-fondamentale}
Let $\mu, \nu : (\Ac,\Dbb) \to \R$ be distributive set functions with $\Dbb$ infinitesimal and $\nu$  positive. 
If  $S \in \Ac$ is a closed and bounded non-empty set and $g$ is the strict derivative of $\mu$ with respect to $\nu$
on $S$,
then   
\begin{equation} \label{mediafor}
\inf_{S} g\leq \frac{\mu (A)}{\nu(A)}\leq \sup_{S} g  
\end{equation}
 for all $A\in \Ac$ with $A\subset S$ and $\nu (A)>0$.\end{theorem}
In the case $\Dbb=\Dbb_{\textrm{int}}$, \textsc{Peano} proves this fundamental theorem  by applying Theorem \ref{teo-F2} to the semi-distributive families 
$\Fo_a:=\{A\in\Ac\, :\,\mu(A)>a\, \nu(A)\}$ and $\Gc_a:=\{A\in\Ac\, :\,\mu(A)<a\, \nu(A)\}$, for real numbers $a$.
Observe that (\ref{mediafor}) amounts to (\ref{Hahn_1})-(\ref{Hahn_2}) and also, indirectly,
to (\ref{hom})-(\ref{sum}).

In \emph{Applicazioni geometriche}, Theorem \ref{teorema-fondamentale} is followed by three corollaries, 
which we summarize into the following:
\begin{corollary} \label{cor-fond} \cite[(1987) p.\,171]{peano87}
Under the same hypothesis as in the previous theorem:
\begin{enumerate}
\item \labelpag{cor-1} if the strict derivative $\frac{d\,\mu}{d\,\nu}$ is a constant $b$ on $S$, then $\mu (A)=b\, \nu(A)$, for all $A\in \Ac$ with 
$A\subset S$;
\item \labelpag{cor-2} if the strict derivative $\frac{d\,\mu}{d\,\nu}$ vanishes at every point of $S$, then $\mu(A)=0$, for all $A\in \Ac$ with 
$A\subset S$; 
\item  \labelpag{due-mis}
if two distributive set functions have equal strict derivatives with respect to $\nu$ on $S$, 
then they are equal on subsets of $S$ belonging to $\Ac$.\,\footnote{
It is evident that properties (\ref{cor-1})-(\ref{due-mis}) are equivalent.
To prove (\ref{due-mis}), \textsc{Peano} shows that the strict derivative of a sum of two distributive set functions
is the sum of their derivatives. 
}
\end{enumerate} 
\end{corollary}

The following fundamental \textsc{Peano}'s result point out the difference of \textsc{Peano}'s approach with respect to both approaches of \textsc{Cauchy} and of \textsc{Lebesgue} (1910).

\begin{theorem}
Under the same hypothesis as in the previous theorem,
if the strict derivative of $\mu$ with respect to $\nu$ exists on $S$, then it is continuous on $S$.
\end{theorem}

The importance of these results is emphasized in \emph{Applicazioni geometriche} by a large amount of evaluations of derivatives of distributive set functions.
As a consequence of the existence of the strict derivative, \textsc{Peano} gives, for the first time, 
several examples of measurable sets.
The most significant examples, observations  and results are listed below.

\begin{enumerate}
\item  \labelpag{hypo-der}
{\it Measurability of the hypograph of a continuous function} \cite[(1887) pp.\,172-174]{peano87}. Let $f$ be a continuous positive real function defined on an interval $[a,b]$, 
 let $\Ac$ be the family of all sub-intervals of $[a,b]$ and 
let $\nu$ be the Euclidean measure on $1$-dimensional intervals. 
Define $\mu_f : \Ac \to \R$ on every $A$ belonging to $\Ac$, by the inner (respectively, the outer) $2$-dimensional measure
(in the sense of Peano-Jordan) of the \emph{positive-hypograph} of $f$, restricted to $A$.\,\footnote{By \emph{positive-hypograph} of $f$ restricted to $A$ we mean the set 
$\{(x,y)\in [a,b]\times R_+ : x\in A$ and $y\leq f(x), \}$, 
where $\R_+ := \{ x\in \R: x\ge 0\}$.} 
In any case, independently of the choice of inner or outer measure, we have
that $\mu_f$ and $\nu$ are distributive set functions with respect to $\Dbb(a,b)$, and that 
 $\frac{\de \mu_f}{\de \nu}(x)=f(x)$ for every $x \in [a,b]$. 
From  (\ref{due-mis}) of Corollary \ref{cor-fond} it follows  that 
the inner measure of the positive-hypograph of the continuous function $f$ coincides with its outer measure;
therefore it is measurable in the sense of Peano-Jordan.

\item 
Analogously, \textsc{Peano} considers continuous functions of two variables and the volume of the \emph{positive-hypograph} \cite[(1887) p.\,175]{peano87}. 

\item \labelpag{area-star}
\emph{Area of a plane star-shaped subset delimited by a continuous closed curve}
\cite[(1887) pp.\,175-176]{peano87}. Consider a 
continuous closed curve that can be described in polar coordinates
in terms of a continuous function 
$\rho:[0,2\pi] \to \R_+$, with $\rho(0)=\rho(2\pi)$. 
Let $\Ac$ be the family of all subintervals of $[0,2\pi]$; and 
for every $A\in\Ac$, let  $\nu(A)$ denote the Euclidean measure  of the area of the circular sector 
$\{(\rho \cos(\theta),\rho\sin(\theta)) : \theta \in A, \rho\in [0,1]\}$. 
Moreover, let  $\mu(A)$ denote  inner (or outer, indifferently) Peano-Jordan $2$-dimensional measure of the set
$\{(\rho\cos(\theta),\rho\sin(\theta)) : \theta \in A, \rho\in [0,\tilde\rho (\theta)]\}$. Then 
the strict derivative
$\frac{\de \mu}{\de \nu}(\theta )$ is equal to $\rho^2 (\theta)$. 
From the fact that this derivative does not depend on the choice of inner or outer measure,
it follows   Peano-Jordan measurability of plane star-shaped sets delimited by continuous closed curves.

\item 
Analogously, \textsc{Peano} considers the volume of a star-shaped set bounded by simple continuous closed surface
\cite[(1887) p.\,177]{peano87}. 

\item \labelpag{cav-pri}{\it Cavalieri's principle between a prism and a spatial figure}
\cite[(1887) pp.\,177-179]{peano87}.
Consider a straight line $r$ in the tri-dimensional space, an unbounded cylinder $P$ parallel to $r$ with polygonal section, and a 
spatial figure $F$. Let $\pi_x$ denote the plane perpendicular to $r$ at the point $x\in r$. Assume Peano-Jordan measurability of all sections of $F$ perpendicular to $r$, namely  
$$
\mu_e(\partial F\cap \pi_x)=0 \quad \text{  for all  } x \in r
\leqno \quad \qquad (*)$$
where $\mu_e$ denotes $2$-dimensional Peano-Jordan outer measure and $\partial F$ denotes the boundary of $F$.
Let $\Ac$ be the family of all segments of $r$. Given a set $A\in\Ac$, 
let $\mu : \Ac \to \R$ denote the outer (or inner, indifferently) $3$-dimensional measure of the set
$\cup_{x\in A}(F\cap \pi _x)$,
and $\nu (A)$ denote  Peano-Jordan $3$-dimensional measure of the set
$\cup_{x\in A}(P\cap \pi_x)$. 
The set functions  $\mu$ and $\nu$ are distributive with respect to the family $\Dbb(r)$ of interval-decompositions of $r$ and
$$\frac{\de \mu}{\de \nu}(x )=\frac{\mu_e(F\cap \pi _x)}{\mu_e(P\cap \pi _x)} \quad \text{for every } x \in r.$$
From the fact that this derivative does not depend on the choice of the inner or outer measure
involved in defining $\mu$, it follows Peano-Jordan measurability of the spatial figure $F$.

\item   {\it Cavalieri's principle between two spatial figures} 
\cite[(1887) p.\,180]{peano87}. 
Consider two spatial figures $F$ and $G$ such that all their sections with planes perpendicular to a given straight line $r$ are Peano-Jordan measurable.
Let $\Ac$ be the family of all segments of $r$.
Given a set $A\in\Ac$, 
let $\mu(A)$ and $\nu(A)$ denote  outer (or inner, indifferently) Peano-Jordan $3$-dimensional measures of the sets
$\cup_{x\in A}(F\cap \pi _x)$ and $\cup_{x\in A}(G\cap \pi _x)$, respectively.
The set functions  $\mu$ and $\nu$ are distributive 
with respect to the family $\Dbb(r)$ of interval-decompositions of $r$ and  
$$\frac{\de \mu}{\de \nu}(x)=\frac{\mu_e(F\cap \pi _x)}{\mu_e(G\cap \pi _x)} \quad \text{for every } x \in r.$$
Hence, from (\ref{cor-1}) it follows the classical Cavalieri's principle: two figures whose 
corresponding sections have equal areas, have the same volume.

\item\labelpag{cav-3} {\it Cavalieri's principle for 3 dimensional figures with respect to one dimensional sections}
\cite[(1887) p.\,180]{peano87}.
Consider a plane $\pi$. 
Let $\Ac $ be the family of all rectangles of $\pi$ and
let $r_x$ be the straight line perpendicular to $\pi$ at $x\in\pi$.  
Moreover, consider a spatial figures $F$ such that for any $x\in \pi$
$$ 
\mu_e(\partial F\cap r_x)=0 \quad \text{ for every } x \in \pi  \leqno \quad \qquad (**) 
$$
where $\mu_e$ denotes the Peano-Jordan $1$-dimensional outer measure and $\partial F$ denotes the
boundary of $F$.
Given a set $Q\in\Ac$, let $\mu (Q)$ denote the outer (or inner, indifferently) measure of the set
$\cup_{x\in Q}(F\cap r_x)$,
and $\nu(Q)$ denote the elementary usual measure of $Q$.
Then $\mu$ and $\nu$ are distributive with respect interval-decompositions of rectangles of $\pi$ and  
$$\frac{\de \mu}{\de \nu}(x )= \mu_e(F\cap r_x) \quad \text{ for every } x \in \pi.$$ 

\item \labelpag{cav-pri-pia}
{\it Cavalieri's principle for 2 dimensional figures}
\cite[(1887) p.\,180]{peano87}.
Analogously to (\ref{cav-pri}), \textsc{Peano} considers Cavalieri's principle for planar figures.

\item 

{\it Derivative of the length of an arc}
\cite[(1887) p.\,181]{peano87}.
In order to compare the length of an arc with the length of its orthogonal projection on a straight line $r$, 
\textsc{Peano} assumes that the orthogonal projection is bijective on a segment $\rho$ of  $ r$,
and that the arc can be parametrized through a function with continuous non null derivative.\,\footnote{The requirement that the derivative of the arc with respect to a parameter
be continuous and non null is expressed by \textsc{Peano} in geometrical terms, namely by requiring
that ``the tangent straight line exists at every point $P$ of the arc, and it is the limit of the 
straight lines passing through two points of the arc, when they tend to $P$''.
\textsc{Peano} was aware that these geometrical conditions are implied by the existence of 
a parametrization with a continuous non-null derivative \cite[(1987) p.\,59, 184]{peano87}.
}
Let $\Ac$ be the family of all closed bounded segments of $\rho$.
For every segment  $A \in \Ac$, let $\mu(A)$ denote the length of the arc whose orthogonal projection over $r$ is  $A$ and let $\nu(A)$ denote the length of $A$. 
Then  
$$\frac{\de \mu}{\de \nu}(x )=\frac{1}{\cos\theta_x}  \leqno \qquad \quad (***) $$
where $\theta_x$ is the angle between $r$ and the straight line that is tangent to the arc at the point 
(of the arc) corresponding to $x \in \rho$.\,\footnote{
Of course, to avoid $\cos\theta_x = 0$ along the arc, \textsc{Peano} assumes that the tangent straight line at every point of the arc is not orthogonal to $r$.
}
\item {\it Derivative of the area of a surface}
\cite[(1887) pp.\,182-184]{peano87}.
By adapting the previous argument, \textsc{Peano} shows that the strict derivative between the area of a surface and its projection on a plane is given by (***), where $\cos \theta$ is the cosinus of the angle between the tangent plane and the projection plane.
\end{enumerate}

\section{Distributive set functions: integral and strict derivative}\label{sez-massdensity}
\textsc{Peano} introduces also the notion of integral with respect to a positive distributive set function.
The \emph{proper integral} of a bounded function $\rho $ on a set $A\in \Ac$ with respect to a positive 
distributive set function $\nu \colon (\Ac, \Dbb) \to \R$, 
is denoted by $\int _A\rho \,\de\nu$ and is defined as the real number such that for any $\Dbb$-decomposition $\{A_i\}_{i=1}^m$
of the set $A$, one has 
$$ \int _A\rho \, \de\nu \ge \rho_1'\nu (A_1)+\rho_2'\nu (A_2)+\dots+\rho_n'\nu (A_m)$$
$$ \int _A\rho \,\de\nu \le \rho_1''\nu (A_1)+\rho_2''\nu (A_2)+\dots+\rho_n''\nu (A_m)$$
where $\rho_1',\rho_2',\dots,\rho_m'$  (respectively $\rho_1'',\rho_2'',\dots,\rho_m''$), are numbers 
defined by
\begin{equation} \label{darboux}
\rho'_i :=\inf_{x\in A_i}\rho(x) \quad\text{and} \quad \rho_i'':=\sup_{x\in A_i}\rho (x),
\end{equation} 
for all $i=1,...,m$.\,\footnote{This clear, simple and general definition of integral with respect to an abstract
positive distributive set function is ignored until the year 1915, when \textsc{Fr\'echet} re-discovers it
in the setting of ``finitely additive'' measures \cite[(1915)]{frechet1915}.
}

\textsc{Peano} defines also the {\it lower} $ \underline\int _{A}\rho \, \de \nu$ 
and the {\it upper} integral 
$ \overline\int _{A}\rho \, \de \nu$
of a bounded function $\rho$ on a set $A\in\Ac$ by
$$ \underline\int _{A}\rho \, \de \nu :=\sup s' \quad \text{  and  } \quad
\overline\int _{A}\rho \, \de \nu :=  \inf s''$$
where 
$s'= \rho_1'\nu (A_1)+\rho_2'\nu (A_2)+\dots+\rho_m'\nu (A_m)$ and $s''=\rho_1''\nu (A_1)+\rho_2''\nu (A_2)+\dots+\rho_m''\nu (A_m)$, where $\rho_i'$ and $\rho_i''$ are defined as in (\ref{darboux})
and $\{A_i\}_{i=1}^m$ runs over $\Dbb$-decompositions of $A$.

In \textsc{Peano}'s terminology, the integrals defined above are called \emph{geometric integrals}.
\textsc{Peano} stresses  the analogy among these integrals and the usual \emph{elementary integral} 
$\int_a^b f(x) \, \de x$ of functions $f$ 
defined over intervals of $\R$.

Using property \eqref{dec_1} of $\Dbb$-decompositions,
\textsc{Peano} shows that the lower integral is always less or equal than the upper integral.
When these values coincide, their common value
is called a proper  integral and is denoted by $\int_A\rho\, \de\nu$. 

Moreover, using properties \eqref{dec_2} and \eqref{dec_3} of $\Dbb$-decompositions, 
\textsc{Peano} shows that the lower integral  $A \mapsto \underline\int _{A}\rho \, \de \nu $ and the  upper integral $A \mapsto \bar\int _{A}\rho \, \de \nu$ are distributive set functions on $\Ac$
with respect to the same family $\Dbb$ of decompositions
\cite[(1887) Theorem\,I, p.\,187]{peano87}. 

In case of $\rho$ continuous, using the property of ``infinitesimality'' of $\Dbb$ (see Definition \ref{def_inf}), \textsc{Peano} shows that the derivative
of both lower and upper integrals with respect to $\nu$ is $\rho$ 
\cite[(1887) Theorem\,II, p.\,189]{peano87}; 
consequently the proper integral $\int_A\rho\, \de\nu$ of a continuous $\rho$ exists
whenever $A$ is closed and bounded 
\cite[(1887) Cor. of Theorem\,II, p.\,189]{peano87}. 

The definitions introduced above allow \textsc{Peano} to realize the mass-density paradigm, i.e., to prove that it is possible to recover a distributive function $\mu$
as the integral of the strict derivative $\frac{d\mu}{d \nu}$  with respect to a positive distributive function $\nu$. \textsc{Peano}'s results can be formulated into the following
\begin{theorem}[\textsc{Peano}'s Theorem on strict derivative of distributive set functions, see {\cite[(1887) Theorem\,14, p.\,171, 
Theorems\,II, III, p.\,189]{peano87}}]\label{peanodev} 
Let $\mu,\nu : (\Ac, \Dbb) \to \R$ be distributive set functions, with $\nu$ positive
and $\Dbb$ infinitesimal. 
Let $S\in\Ac$ be a closed bounded set and $\rho:S\to\R$ a function. 
The following properties are equivalent:
\begin{enumerate}
\item $\rho$ is the strict derivative $\frac{d\,\mu}{d\,\nu}$ of $\mu$ with respect to $\nu$ on $S$;
\item $\rho$ is continuous and $\mu (A)=\int_A\rho \, \de \nu$
%\footnote{If $\rho$ is continuous at a point $x$ then the derivative of the function 
%$A\mapsto \int_A\rho \de\nu$ is equal to $\rho$} 
for any $A\subset S$, $A\in\Ac$.
\end{enumerate}
\end{theorem}

\textsc{Peano} applies Theorem \ref{peanodev} to the list of examples of strict derivatives 
of distributive set functions of \S\,\ref{sez-derivata} and obtains the following results.
\begin{enumerate}
\item \emph{Fundamental theorem of integral calculus for continuous functions} \cite[(1887) pp.\,191-193]{peano87}.
Consider a continuous function $f$ on $\R$ and let $F$ be a primitive of $f$.
Define $\mu$ and $\nu$ over the family $\Ac$ of closed bounded intervals $[a,b]$ of $\R$ by
$\mu([a,b]):=F(b) - F(a)$ and $\nu([a,b]):=b-a$. 
Observe that both $\mu$ and $\nu$ are distributive set functions with respect to $\Dbb(\R)$ and
$$\frac{\de\mu}{\de\nu}(x) = 
\lim_{\begin{subarray}{c}a,b\to x \\ a\neq b\end{subarray}}
\frac{F(b) - F(a)}{b-a} = f(x) $$
since $F$ has continuous derivative.\,\footnote{
\textsc{Peano} observes that continuity of derivative of $F$ is a necessary and sufficient condition
to have the existence of $\frac{\de\mu}{\de\nu}$.}
Therefore, by Theorem \ref{peanodev}, \textsc{Peano} obtains
$$F(b)-F(a)=\mu([a,b]) =\int_{[a,b]} f\, d\, \nu = \int_a^b f(x)\, \de x  \,.$$

\item {\it Calculus of an integral as a planar area} \cite[(1887) pp.\,193-195]{peano87}. 
The elementary integral of a continuous positive function is 
Peano-Jordan measure of the positive hypograph of the function. This is an immediate
application of Theorem \ref{peanodev} to the setting (\ref{hypo-der}).
\item {\it Cavalieri's formula for planar figures} \cite[(1887) p.\,195]{peano87}. 
Let us suppose that $C\subset\R^2$, $C_x :=\{y\in\R\, :\; (x,y)\in C\}$ and
$(\partial C)_x  := \{y\in\R\, :\; (x,y)\in \partial C\}$ 
for every $x \in \R$. 
Assume that for any $x$ the set $(\partial C)_x$ has vanishing outer measure. 
As a consequence of Theorem \ref{peanodev} and the two-dimensional version of \emph{Cavalieri's principle}
(\ref{cav-pri-pia}) (see \cite[(1887) p.\,180]{peano87}), it follows that
the measure of the part of the figure $C$, bounded by the abscissas $a$ and $b$, is equal to  
$$\int_a^b \mu_e (C_x) \, \de x $$
where $\mu_e$ denotes  outer Peano-Jordan one-dimensional measure.

\item {\it Area of a plane star-shaped subset delimited by a continuous closed curve} 
\cite[(1887) p.\,199]{peano87}.
In the setting of example (\ref{area-star}), \textsc{Peano} shows that the area 
of the sector between the angles $\theta_0$ and $\theta_1$, delimited by a
curve described in polar coordinates by $\rho$, is equal to 
$$\frac{1}{2}\int_{\theta_0}^{\theta_1} \rho(\theta)^2\de \theta \,.$$

\item {\it Cavalieri's formula for volumes} \cite[(1887) p.\,221]{peano87}.
In the setting (\ref{cav-3}), let's define $F_x := \{(y,z) \in\R^2 : (x,y,z)\in F\}$ and
$(\partial F)_x := \{(y,z)\in\R^2 : (x,y,z)\in \partial F\}$.
Assume that for any $x$, 
the set $(\partial F)_x$ has vanishing outer measure. From Theorem \ref{peanodev}, \textsc{Peano} shows that
the volume of the part of the figure $F$,
delimited by the planes $x=a$ and $x=b$, is equal to  $$\int_a^b  \mu_e(F_x)\,\de  x \, $$
where $\mu_e$ denotes  outer  Peano-Jordan two-dimensional measure.

\end{enumerate}

\section{Coexistent magnitudes in Lebesgue and Peano's derivative}\label{sez-comments}

\textsc{Lebesgue} gives a final pedagogical
\footnote{\textsc{Lebesgue} says in \cite[(1931) p.\,174]{lebesgue1931}:
\begin{quote}
[\dots] depuis trente ans [d' enseignement] [\dots] on ne s'\'etonnera pas que l'id\'ee
me soit venue d'\'ecrire des articles de nature p\'edagogique; si j'ose employer ce qualificatif
que suffit ordinairement pour faire fuir les math\'ematiciens.
\translation{%
[\dots] in the thirty years [of teaching] [\dots] it is not at all surprising that the idea 
should occur to me of writing articles on a pedagogical vein; if I may use an expression
which usually puts mathematicians to flight. (transl. \textsc{May} \cite[(1966) p.\,10]{lebesgue_may})}
\end{quote}
}
 exposition of his measure theory
in \emph{La mesure des grandeurs} \cite[(1935) p.\,176]{lebesgue1935}, 
by referring directly to \textsc{Cauchy}'s \emph{Coexistent magnitudes}:\,\footnote{The five parts of the essay \emph{La mesure des grandeurs} have been published in \emph{L'Enseignement math\`ematique} during the years 1931-1935.
An english translation \emph{Measure and the Integral} of \emph{La mesure des grandeurs} is due to Kenneth O.~May \cite[(1966)]{lebesgue_may}.}
\begin{quote}
La th\'eorie des grandeurs qui constitue le pr\'ec\'edent 
chapitre avait \'et\'e pr\'epar\'ee par des recherches de Cauchy, sur 
ce qu'il appelait des grandeurs concomitantes [\emph{sic}], par les travaux
destin\'es \`a \'eclaircir les notions d'aire, de volume, 
de mesure [\dots].\,\footnote{\translation{%
The theory of magnitudes forming the subject of the preceding chapter was
prepared by researches of Cauchy on what he called concomitant magnitudes,
by studies destined to clarify the concepts of area, volume, and measure [\dots]
(transl. \textsc{May} \cite[(1966) p.\,138]{lebesgue_may})}}
\end{quote}

\textsc{Lebesgue} is aware of the obscurity of the concepts that are present in  \textsc{Cauchy}'s 
 \emph{Coexistent magnitudes}, starting by the meaning of the term \emph{magnitude} itself.
In this respect, in order to put on a solid ground the ideas of \textsc{Cauchy},
\textsc{Lebesgue} was compelled to pursuit an approach similar to that \textsc{Peano}: in fact 
he defines a ``magnitude'' as a set function on a family of sets $\Ac$, requires infinitesimality of $\Ac$ (in the sense that 
every element of $\Ac$ can be \emph{r\'eduit \`a un point par diminutions successives}),
and   additivity properties that he express
in \emph{La mesure des grandeurs} \cite[(1934) p.\,275]{lebesgue1934} in these words:
\begin{quote}
Si l'on divise un corps $C$ en un certain nombre de corps
partiels $C_1, C_2,$ $ \dots, C_p$, et si la grandeur $G$ est, pour ces corps, $g$ d'une
part, $g_1, g_2, \dots, g_p$ d'autre part, on doit avoir: $g=g_1 + g_2 + \dots + g_p$.\,\footnote{\translation{%
If a body $C$ is partitioned into a certain number of sub-bodies $C_1, C_2, \dots, C_p$ and 
if for these bodies the magnitude $G$ is $g$ on the one hand and 
$g_1, g_2, \dots, g_p$ on the other, we must have $g= g_1 + g_2 + \cdots + g_p$.
(transl. \textsc{May} \cite[(1966) p.\,129]{lebesgue_may})}

\textsc{Lebesgue} observes that in order to make this condition rigorous, it would be necessary
to give a precise meaning to the words \emph{corp} and \emph{partage de la figure totale
en parties} \cite[(1934) p.\,275-276]{lebesgue1934}. 
Moreover he observes that \emph{diviser un corps} may be interpreted in
different ways \cite[(1934) p.\,279]{lebesgue1934}.}
\end{quote}

In \emph{La mesure des grandeurs} \textsc{Lebesgue} considers the operations of integration and differentiation by 
presenting these topics in a new form with
respect to his fundamental and celebrated paper \emph{L'int{\'e}gration des fonctions discontinues} 
\cite[(1910)]{lebesgue1910}.

\textsc{Lebesgue} theory of differentiation of 1910 concerns absolutely continuous $\sigma$-additive measures
on Lebesgue measurable sets. On the contrary, twenty-five years later
in \emph{La mesure des grandeurs} of 1935
\begin{itemize}
\item
$\sigma$-additive set functions are replaced by 
continuous \footnote{It is not easy to give in a few words a definition of the concept of continuity according to \textsc{Lebesgue}:
such a continuity is based on a convergence of sequences of sets that in the relevant cases
coincides with the convergence in the sense of Hausdorff. 
We recall that a sequence of sets $\Delta_n$ \emph{converges to $\Delta$ in the sense of Hausdorff} if for all $\epsilon >0$ 
there exists $n_0$ such that $\Delta_n\subset B_\epsilon (\Delta)$ and $\Delta\subset B_\epsilon (\Delta_n)$ for all $n>n_0$, where $B_\epsilon (A):= \{x \in \R^n : \text{ there exists } a \in A 
\text{ such that } \|x-a\|<\epsilon \}$. 
Therefore, a set function $f$ is said to be \emph{continuous} if for any
$\Delta_n$ and $\Delta$ Peano-Jordan measurable sets, we have that
$\lim_{n\to \infty} f (\Delta_n) = f(\Delta)$, whenever $\Delta_n$ converges to $\Delta$
in Hausdorff sense. 
} 
additive
\footnote{\textsc{Lebesgue} writes in \cite[(1935) p.\,185]{lebesgue1935}:
\begin{quote}
[\dots] nous supposerons cette fonction [$f$] \emph{additive}, c'est-a-dire telle que,
si l'on divise $\Delta$ en deux domaines quarrables $\Delta_1$ et $\Delta_2$
on ait $f(\Delta) =  f(\Delta_1) + f(\Delta_2)$.
\translation{%
[\dots] let us assume that this function is \emph{additive}; that is, it is such that, if we partition
$\Delta$ into two quadrable domains $\Delta_1$ and $\Delta_2$, we have
$f(\Delta) = f(\Delta_1) + f(\Delta_2)$. (transl. \textsc{May} \cite[(1966) p.\,146]{lebesgue_may})}
\end{quote}
}
 measures;
\item
absolutely continuous measures become set functions with 
bounded-derivati\-ve~\footnote{A set function $f$ has a \emph{bounded-derivative} 
with respect to Peano-Jordan $n$-dimensional measure ${\vol}_n$ if there exists a constant $M$ such that $|f(\Delta)|\leq M \, {\vol}_n(\Delta)$ for any Peano-Jordan measurable set $\Delta$.
A set function with bounded-derivatives is called \emph{uniformly Lipschitzian} by \textsc{Picone}
\cite[(1923) vol.\,2, p.\,467]{picone}.
} (\emph{\`a nombres d\'eriv\'es born\'es});
\item
Lebesgue measurable sets are replaced by Jordan-Peano measurable subsets of a given bounded set.
\end{itemize}

Let $K$ be a bounded closed subset of  Euclidean space $\R^n$, let $\Ac_K$ be the family of 
Jordan-Peano measurable (\emph{quarrables}) subsets of $K$ and let $V$ be a 
positive, continuous, additive
set function on $\Ac_K$ with bounded-derivative. 
Then \textsc{Lebesgue} introduces a definition of derivative.
The \emph{uniform-derivative} (\emph{d\'eriv\'ee \a convergence uniforme}) $\varphi$ of a set function $f$ with respect to $V$, is defined as the function $\varphi:K\to \R$ such that, for every $\epsilon>0$, there exists  $\eta >0$ such that
\begin{equation}
\Big|\frac{f(\Delta)}{V(\Delta)}-\varphi (x)\Big|<\epsilon 
\end{equation} 
for all $x\in K$ and $\Delta \in \Ac_K$ with $x\in \Delta \subset B_\eta (x)$. 
It is clear that \textsc{Lebesgue}'s new notion of uniform-derivative is strictly related to \textsc{Peano}'s one.
In fact, \textsc{Lebesgue} observes that the uniform-derivative is continuous whenever it exists;
moreover, he defines the integral 
\begin{equation}
\int_K \varphi \,\de V
\end{equation}
of a \emph{continuous function} $\varphi$ with respect to $V$. His definition of integral 
\cite[(1935) pp.\,188-191]{lebesgue1935} is rather intricate with respect to that of \textsc{Peano}.

It is worthwhile noticing that \textsc{Lebesgue} recognizes the relevance of the notion of an integral with respect to  set functions.
\textsc{Lebesgue}, not acquainted with previous \textsc{Peano}'s contributions, assigns the priority of this notion to \textsc{Radon} \cite[(1913)]{radon1913}.
On the other hand, \textsc{Lebesgue} notices that the integral with respect to set functions
was already present in Physics
\footnote{\textsc{Lebesgue} gives several examples of this. For instance, the evaluation of the heath quantity necessary to increase the temperature of a body as integral of the specific heath with respect to the mass.}
and express his great surprise in recovering in \textsc{Stieltjes}'s integral
\cite[(1894)]{stieltjes} an instance of
integral with respect to set functions; Lebesgue writes
\cite[(1926) p.\,69-70]{lebesgue1926}:
\begin{quote}
Mais son premier inventeur, Stieltj\`es, y avait \'et\'e conduit par des 
recher\-ches d'analyse et d'arithm\'etique et il l'avait pr\'esent\'ee sous une forme
purement analytique qui masquait sa signification physique;
s\`\i\  bien qu'il a fallu beaucoup d'efforts pour comprendre et 
conna\^itre ce qui est maintenant \'evident. L'historique de ces
efforts citerait les nom de F.~Riesz, H.~Lebesgue, W.H.~Young, M.~Fr\'echet,
C. de la Vall\'e-Poussin; il montrerait que nous avons rivalis\'e en ing\'eniosit\'e, en perspicacit\'e, mais aussi en aveuglement.\,\footnote{\translation{%
But its original inventor, Stieltjes, was led to
it by researches in analysis and theory of number and he presented it in a purely
analytical form which masked its physical significance, so much so that it required 
a much effort to understand and recognizes what is nowadays obvious. The history
of these efforts includes the works of
F.~Riesz, H.~Lebesgue, W.H.~Young, M.~Fr\'echet,
C. de la Vall\'e-Poussin. It shows that we were rivals in ingenuity, in insight, but also
in blindness.
(transl. \textsc{May} \cite[(1966) p.\,190]{lebesgue_may})
}}
\end{quote}

The first important theorem presented by \textsc{Lebesgue} is the following
\begin{theorem} 
Let $K$ be a bounded closed subset of $\R^n$, $\varphi:K\to \R$ a continuous function and $V$ a positive additive continuous set function with bounded-derivative. 
Then the integral $\Delta \mapsto \int_\Delta \varphi \,\de V$ with $\Delta \in \Ac$ is the unique additive set function with bounded-derivative which has $\varphi$ as uniform-derivative 
with respect to $V$.\,\footnote{The proof is rather lengthy, as \textsc{Lebesgue} included in  it the definition of integral as well as the theorem of average value.}
\end{theorem}

The main applications of this theorem, given by \textsc{Lebesgue} in \emph{La mesure des grandeurs}
\cite[(1935) p.\,176]{lebesgue1935}, concern:
\begin{enumerate}
\item the proof that multiple integrals can be given in terms of simple integrals; 
\item  the formula of change of variables;\,\footnote{Lebesgue uses the implicit function theorem.}
\item  several formulae for oriented integrals (Green's formula, length of curves and area of surfaces).
\end{enumerate}

The uniform-derivative defined by \textsc{Lebesgue} is, as observed above, a continuous function, and coincides exactly with \textsc{Peano}'s strict derivative.  Through a different and more difficult path \footnote{The exposition of 1935 is elementary, but more lengthy and difficult than those presented by \textsc{Lebesgue} in 1910.
Surprisingly, the  terms \emph{domain, decomposition, limit, additive, continuous} are used by \textsc{Lebesgue} in a supple way. } than \textsc{Peano}'s one, \textsc{Lebesgue} rediscovers the importance of the continuity of the derivative. 
In \textsc{Lebesgue}'s works there are no references to the contributions of \textsc{Peano}
concerning differentiation of set functions. 

Several years before \emph{La mesure des grandeurs} of 1935, \textsc{Lebesgue} in \cite[(1926)]{lebesgue1926} outlines his contribution to the notion
of integral. In the same paper he mentions 
\textsc{Cauchy}'s \emph{Coexistent magnitudes} in the setting of derivative of measures.
Moreover he cites \textsc{Fubini}'s and \textsc{Vitali}'s works of 1915 and 1916 (published 
by Academies of Turin and of Lincei) in the context of the general problem of primitive functions.

More precisely, in 1915, the year of publication of \textsc{Peano}'s paper \emph{Le grandezze coesistenti} \cite{peano1915},
\textsc{Fubini} \cite[(1915)]{fubini1915b,fubini1915a} and \textsc{Vitali}  
\cite[(1915, 1916)]{vitali1915,vitali1916} introduce a definition of  derivative of ``finitely additive measures''
\footnote{\textsc{Fubini}'s first paper \cite{fubini1915a} is presented by C.\,\textsc{Segre}  at the Academy of Sciences of Turin on January 10, 1915. 
In the same session, \textsc{Peano}, Member of the Academy, presents a multilingual dictionary and a paper written by one of his students,  \textsc{Vacca}.
\textsc{Segre}, on April 11, 1915, presents, as a Member, a second paper of \textsc{Fubini} \cite{fubini1915b} to {\it Accademia dei Lincei}.
In the session of the Academy of Turin of June 13, 1915, \textsc{Peano} presents his paper \emph{Le grandezze coesistenti}. 
Moreover \textsc{Segre} presents two papers by \textsc{Vitali} \cite[(1915)]{vitali1915} and
\cite[(1916)]{vitali1916}
to Academy of Turin on November 28, 1915 and to Academy of Lincei on May 21, 1916, respectively.  

There is a rich correspondence between \textsc{Vitali} and \textsc{Fubini}. In the period March-May 1916
\textsc{Fubini} sends three letters to \textsc{Vitali} 
(transcribed in \emph{Selected papers} of \textsc{Vitali} \cite[pp.\,519-520]{vitali-opere}), concerning
differentiation of finitely additive measures and related theorems. In particular \textsc{Fubini} suggests 
\textsc{Vitali} to quote \textsc{Peano}'s paper \cite[(1915)]{peano1915} and to compare alternative definitions
of derivative.
In \emph{Selected papers} of \textsc{Vitali} it is also possible to find six letters by \textsc{Peano} to \textsc{Vitali}.
Among them, there is letter of March 21, 1916 concerning \textsc{Cauchy}'s coexistent magnitudes; \textsc{Peano} writes:
\begin{quote}
Grazie della sua nota \cite[(1915)]{vitali1915}. Mi pare che la dimostrazione che Ella
d\`a, sia proprio quella di Cauchy, come fu rimodernata da G. Cantor, e poi da me, e di cui
trattasi nel mio articolo, Le grandezze coesistenti di Cauchy, giugno 1915, e di cui debbo avere inviato
copia.  
\end{quote}

[\![Thanks for your paper \cite[(1915)]{vitali1915}. In my opinion your proof coincides with the one
given by Cauchy, as formulated by Cantor and by myself in my paper 
``Coexistent magnitudes of Cauchy'' (June 1915), that I sent you.]\!]

To our knowledge, \textsc{Fubini} \cite[(1915)]{fubini1915b,fubini1915a} and \textsc{Vitali}
\cite[(1915, 1916)]{vitali1915,vitali1916} are not cited by other authors, with the exception of 
\textsc{Banach} \cite[(1924) p.\,186]{banach1924}, who refers to \textsc{Fubini}  \cite[(1915)]{fubini1915b}.
}, oscillating themselves 
between  definitions \emph{\`a la Cauchy}  and \emph{\`a la Peano}. 

\textsc{Vitali}, in his second paper \cite{vitali1916}, refers to the \emph{Coexistent magnitudes} of \textsc{Cauchy}, and 
presents a comparison among the notions of derivative given by \textsc{Fubini}, himself, \textsc{Peano} and the one of \textsc{Lebesgue} of 1910, emphasizing the continuity of the \textsc{Peano}'s strict derivative. \textsc{Vitali} writes in 
\cite[(1916)]{vitali1916}:
\begin{quote}
Il Prof.\,G.\,Peano nella Nota citata [\emph{Le grandezze coesistenti}] e in un'altra sua
pubblicazione anteriore [\emph{Applicazioni geometriche}], si occupa dei teoremi di Rolle e della 
media e ne indica la semplice dimostrazione nel caso in cui la derivata [della funzione di insieme $f$]
in $P$ sia intesa come il limite del rapporto  di $\frac{f(\tau)}{\tau}$, dove $\tau$ \`e un campo
qualunque che pu\`o anche non contenere il punto $P$.

L'esistenza di tale simile derivata finita in ogni punto porta difatti la continuit\`a
[della derivata medesima].\,\footnote{\translation{%
Prof.\,Peano, in the cited Paper [\emph{Le grandezze coesistenti}] and in a previous publication 
[\emph{Applicazioni geometriche}] deals with  Rolle's and mean value theorems,
pointing out a simple proof, valid in the case in which the derivative  [of the set function $f$],  in
a given point $P$,  is the limit of the ratio $\frac{f(\tau)}{\tau}$, where $\tau$ is a set
that might not contain the point $P$.}}
\end{quote}

This proves that since 1926 \textsc{Lebesgue} should have  been aware of
\textsc{Peano}'s derivative and of its continuity.\,\footnote{We can ask how much \textsc{Lebesgue} was aware of the contributions of \textsc{Peano}.
In many historical papers the comment of \textsc{Kennedy} \cite[(1980) p.\,174]{peano_vita}, a well known 
biographer of \textsc{Peano}, occurs:
\begin{quote}
Lebesgue acknowledged Peano's influence on his own development.
\end{quote} 
In our opinion \textsc{Peano}'s influence on \textsc{Lebesgue} is  relevant but  sporadic.
After a reading of \textsc{Lebesgue}'s works, we have got the feeling that his knowledge of
\textsc{Peano}'s contributions was restricted to two papers on the definition of area 
and on Peano's curve.
}

Undoubtably, the contributions of \textsc{Peano} and \textsc{Lebesgue} have
a pedagogical and mathematical relevance in formulating a definition of derivative
having the property of continuity whenever it exists. Surprisingly these contributions are not known.

Rarely the notion of derivative of set functions is presented and used in
educational texts. 

An example is provided by \emph{Lezioni di analisi matematica} 
of \textsc{Fubini}. There are several editions of these \emph{Lezioni}: starting by the 
second edition \cite[(1915)]{fubiniB1915}, \textsc{Fubini} introduces a derivative 
\emph{\`a la Peano} of additive set functions in order to build a basis for integral calculus in one or several variables. Nevertheless, in his \emph{Lezioni}, \textsc{Fubini} assumes
continuity of its derivative  as an additional property. Ironically, \textsc{Fubini} is aware of continuity of \textsc{Peano}'s derivative, whenever it exists; this is clear from two letters of 1916 that he sent to \textsc{Vitali} \cite[p.\,518-520]{vitali-opere}; in particular, in the second letter, about the \textsc{Peano}'s paper \emph{Grandezze coesistenti} 
\cite[(1915)]{peano1915}, he writes:
\begin{quote}
Sarebbe bene citare [l'articolo di] Peano e dire che, se la derivata esiste e per
calcolarla in [un punto] $A$ si adottano anche dominii che tendono ad
$A$, pur non contenendo $A$ all'interno, allora la derivata \`e continua.\,\footnote{[\![It would be important to cite the paper of Peano, saying that, whenever the 
derivative exists and its evaluation is performed by considering domains that approach $A$, without requiring that the point $A$ belongs to the domains themselves, 
then the derivative is continuous.%
]\!]}
\end{quote}

The notion of derivative of set function is also exposed in the
textbooks \emph{Lezioni di analisi infinitesimale} of \textsc{Picone} \cite[(1923) vol.\,II, p.\,465--506]{picone},
in \emph{Lezioni di analisi matematica} of \textsc{Zwirner} \cite [(1969), pp.\,327-335]{zwirner}
and in \emph{Advanced Calculus} of R.\,C. and E.\,F. \textsc{Buck} \cite[(1965)]{buck}.
In the book of \textsc{Picone}, a definition of derivative 
\emph{\`a la Cauchy} of ``additive'' set functions is given;\,\footnote{Significant instances of additive set functions in the sense of \textsc{Picone}
are outer measure of Peano-Jordan on all subsets of $\R^n$   and lower/upper integrals of functions
with respect to arbitrary domain of integration \cite[(1923) vol.\, II, p.\,356-357, 370-371]{picone}.
The family of decompositions that leads to the notion of additive set function in the sense of \textsc{Picone} is clearly defined on page 356-357 of his book \cite{picone} and includes
the family of decompositions (\ref{dec-picone}) and (\ref{dec-picone2}).
} 
it represents an improvement of \textsc{Cauchy}, \textsc{Fubini} and \textsc{Vitali} definitions. 
Of course, his derivative is not necessarily a continuous function.  Whenever the derivative is continuous, 
\textsc{Picone} states a fundamental theorem of calculus, and applies it to the change of
variables in multiple integrals.
In the book of \textsc{Zwirner} the notion of  derivative \emph{\`a la Peano} of set functions is introduced, without mentioning \textsc{Peano} and, unfortunately, without providing any application. In the third book, R.\,C. and E.\,F. \textsc{Buck} introduce in a clear way a simplified notion of the uniform-derivative of \textsc{Lebesgue} (without mentioning him), and they apply it to obtain the basic formula for
the change of variables in multiple integrals.

\section{Appendix}

All articles of \textsc{Peano} are collected in \emph{Opera omnia} \cite{peano_omnia}, a CD-ROM edited by C. S. \textsc{Roero}. Selected works of \textsc{Peano} were assembled and commented in \emph{Opere scelte} \cite{peano_opere} by \textsc{Cassina}, a pupil of \textsc{Peano}. For a few works there are  English translations in \emph{Selected Works} \cite{peano_english}. Regrettably, fewer \textsc{Peano}'s papers have a  public URL and are freely downloadable.

For reader's convenience, we provide a chronological list of some
mathematicians mentioned in the paper, together with biographical sources.

\texttt{Html} files  with biographies of mathematicians listed below with an asterisk can be attained at University of St
Andrews's web-page 
\begin{center}\texttt{http://www-history.mcs.st-and.ac.uk/history/\{Name\}.html}
\end{center}

\textsc{Kepler}, Johannes (1571-1630)*

\textsc{Cavalieri}, Bonaventura (1598-1647)*

\textsc{Newton}, Isaac (1643-1727)*

\textsc{Mascheroni}, Lorenzo (1750-1800)*

\textsc{Cauchy}, Augustin L. (1789-1857)*

\textsc{Lobachevsky}, Nikolai I. (1792-1856)*

\textsc{Moigno} Fran\c cois N.~M. (1804-1884),
see \emph{Enc. Italiana}, Treccani, Roma, 1934

\textsc{Grassmann}, Hermann (1809-1877)*

\textsc{Serret}, Joseph A. (1819-1885)*

\textsc{Riemann}, Bernhard (1826-1866)*

\textsc{Jordan}, Camille (1838-1922)*

\textsc{Darboux}, Gaston (1842-1917)*

\textsc{Stolz}, Otto (1842-1905)*

\textsc{Schwarz}, Hermann A. (1843-1921)*

\textsc{Cantor}, Georg (1845-1918)*

\textsc{Tannery}, Jules (1848-1910)*

\textsc{Harnack}, Carl (1851-1888),  see May \cite[(1973) p.\,186]{may1973}

\textsc{Stieltjes}, Thomas J. (1856-1894)*

\textsc{Peano}, Giuseppe (1858-1932)*, see  \cite{peano_vita}

\textsc{Young}, William H. (1863-1942)*

\textsc{Segre}, Corrado (1863-1924)*

\textsc{Vall\'ee Poussin} (de la), Charles (1866-1962)* 

\textsc{Hausdorff}, Felix (1868-1942)*

\textsc{Borel}, Emile (1871-1956)*

\textsc{Vacca}, Giovanni (1872-1953)*

\textsc{Carath\'{e}odory}, Constantin (1873-1950)*

\textsc{Lebesgue}, Henri (1875-1941)*

\textsc{Vitali}, Giuseppe (1875-1932)*

\textsc{Fr\'{e}chet}, Maurice (1878-1973)*

\textsc{Fubini}, Guido (1879-1943)*

\textsc{Riesz}, Frigyes (1880-1956)*

\textsc{Tonelli}, Leonida (1885-1946)*

\textsc{Picone}, Mauro (1885-1977), see \texttt{http://web.math.unifi.it}

\textsc{Ascoli}, Guido (1887-1957),
see May \cite[(1973) p.63]{may1973}

\textsc{Radon}, Johann (1887-1956)*

\textsc{Nikodym}, Otton (1887-1974)*

\textsc{Bouligand}, George (1889-1979),
\emph{see} \texttt{http://catalogue.bnf.fr}

\textsc{Banach}, Stefan (1892-1945)*

\textsc{Kuratowski}, Kazimierz (1896-1980)*

\textsc{Cassina}, Ugo (1897-1964), 
see  Kennedy \cite[(1980)]{peano_vita}

\textsc{Cartan}, Henri (1904-2008)*

\textsc{Dieudonn\'e}, Jean A.~E. (1906-1992)*

\textsc{Choquet}, Gustave (1915-2006),
see \emph{Gazette des Math.} v111:74-76, 2007

\textsc{May} Kenneth O. (1915-1977), see \cite[p.\,479]{dauben_scriba}

\textsc{Medvedev} F\"edor A. (1923-1993), see \cite[p.\,482]{dauben_scriba}

\bibliographystyle{plain}

\end{document}